\documentclass[12pt]{article}
\usepackage{color}
\usepackage{verbatim}
\usepackage[polutonikogreek,english]{babel}
\usepackage{booktabs,caption}
\usepackage{multicol,multirow,tabularx}
\usepackage{siunitx}
\usepackage{relsize}
\usepackage{multirow}
\usepackage{hyperref}
\usepackage{float}
\usepackage{lineno}
\usepackage[left=2cm, right=5cm, top=2cm]{geometry}
\usepackage{eurosym}
\usepackage{graphics}
\usepackage{epsfig}
\usepackage{amssymb}
\usepackage{amsmath}
\usepackage{adjustbox}
\usepackage{xcolor}
\usepackage{amsthm}
\usepackage{enumitem}

\newtheoremstyle{break}
  {\topsep}{\topsep}%
  {\itshape}{}%
  {\bfseries}{}%
  {\newline}{}%
\theoremstyle{break}
\newtheorem{theorem}{Theorem}[section]
\newtheorem{proposition}{Proposition}[section]
\newtheorem{corollary}{Corollary}[section]
\newtheorem{lemma}{Lemma}[section]
\newtheorem{remark}{Remark}[section]

\usepackage[square,sort,comma,numbers,compress]{natbib}

\DeclareCaptionFont{ninept}{\fontsize{9pt}{11pt}\selectfont #1}
\usepackage[left=2cm, right=5cm, top=2cm]{geometry}
\begin{document}

\title{Global phase-plane organization \\
of a toxin-mediated plant-herbivore model}
\author{Grif\'o Gabriele\\
{\footnotesize Department of Mathematical, Computer, Physical and Earth Sciences, University of Messina,} \\
{\footnotesize V.le F. Stagno D'Alcontres 31, I-98166 Messina, Italy.}\\
{\small gabgrifo@unime.it}
}
\date{}
\maketitle

\begin{abstract}

Plant-herbivore interactions with chemical plant defence mechanisms can generate complex population dynamics through the nonlinear interplay between resource availability, herbivore pressure, and defence-induced feedbacks. In this work, we investigate the dynamical organization of a plant-herbivore model characterized by a non-monotone functional response accounting for inhibitory effects on grazing activity. In particular, we perform a qualitative and geometric analysis of the associated homogeneous system, characterizing the dynamical regimes supported by the model and the transitions between them. The results reveal the emergence of oscillatory behaviour through a Hopf bifurcation and the occurrence of stable periodic solutions. We further characterize the global phase-plane organization, providing a complete classification of equilibrium configurations and heteroclinic structures supported by the model. In addition, we identify the conditions under which homoclinic connections may occur and discuss their role in shaping the global dynamics. The analysis provides a unified description of the phase-plane organization and highlights how nonlinear defence feedbacks generate qualitatively distinct ecological outcomes, including extinction, stable coexistence, persistent oscillations, and global bifurcation scenarios organized by heteroclinic and homoclinic structures.\\

\textbf{Keywords}: plant-herbivore model; heteroclinic and homoclinic connections; Hopf bifurcation; population oscillations.\\

\end{abstract}

\section{Introduction}
\label{sec:intro}

Plant-herbivore interactions are shaped by the interplay between resource growth, grazing pressure, and defence mechanisms developed by plants to reduce herbivore exploitation. Such interactions may generate a wide variety of ecological outcomes, including species persistence, extinction scenarios, bistability, and sustained population oscillations \cite{Rosenzweig1963, May1972, Gilpin1973, Crawley1983, Krupa1990, Karban1997, Hofbauer1998, Fussmann2000, Murdoch2003, Stiling2007, Castellanos2019}. 

Among the various defence strategies adopted by plants, the production of toxic secondary metabolites plays a particularly important role, as it may directly affect herbivore feeding activity and alter the balance between resource availability and consumption pressure. Therefore, the ecological consequences of toxin-mediated interactions have been investigated through a series of mathematical models incorporating defence-dependent functional responses. Indeed, previous studies have shown that plant toxicity may substantially modify the existence and stability of coexistence states, generate complex bifurcation structures, and promote oscillatory population dynamics \cite{Li2006, Feng2008, Liu2008, Castillo2012, Feng2013, Xiang2021}. More recently, attention has also extended to spatially distributed formulations, where diffusion and cross-diffusion mechanisms have been shown to induce pattern formation and spatiotemporal self-organization \cite{Xiang2021,Grifo2027}. 

Despite these advances, most studies have focused on the existence and stability of equilibria, local bifurcations, and numerical exploration of the resulting dynamics. Although several important global results are available \cite{Liu2008, Castillo2012, Feng2013}, the geometric mechanisms organizing the phase portrait remain less explicitly characterized. In particular, the role played by invariant manifolds and by heteroclinic and homoclinic structures in shaping the long-term dynamics of toxin-mediated plant-herbivore systems has not yet been systematically investigated. A geometric description of these invariant structures is essential for understanding how local bifurcation phenomena are embedded into the global organization of the phase plane and for clarifying the dynamical scenarios observed in numerical studies.

Motivated by these considerations, in this work we revisit a toxin-mediated plant-herbivore model characterized by a non-monotone functional response accounting for inhibitory effects on grazing activity. The novelty of the present work is twofold. First, we provide a complete geometric classification of the global phase portraits arising in the model. Second, we perform a detailed analysis of the parameter plane, showing how changes in the model parameters reorganize the global geometry of the system. In particular, we characterize the bifurcation scenarios connecting different dynamical regimes, including the onset of periodic solutions through Hopf bifurcation and the appearance of heteroclinic and homoclinic structures \cite{Aguirre2014, Wei2017, Shang2022}. This unified geometric framework reveals how the interaction between plant toxicity and herbivore exploitation shapes the global dynamics, leading to distinct ecological outcomes ranging from extinction to sustained oscillatory behavior.

The manuscript is organized as follows. In Section \ref{sec:model} we introduce the model and discuss its biological assumptions. Section~\ref{sec:equilibria} is devoted to the analysis of equilibria, local stability properties, and Hopf bifurcation. The global phase-plane organization, including heteroclinic and homoclinic structures, is investigated in Section~\ref{sec:global}. Finally, conclusions and perspectives for future research are presented in Section~\ref{sec:disc}.

\section{Model formulation and geometric preliminaries}
\label{sec:model}

We consider the toxin-mediated plant--herbivore model proposed in \cite{Li2006,Liu2008,Feng2008, Grifo2027}. The dimensionless equations are 

\begin{equation} 
\label{eq:model} 
\begin{aligned} 
\dot B &= B\left(1-\frac{B}{\alpha}\right) -mC(B)H, \\[1mm] 
\dot H &= \left(mC(B)-\theta\right) H, 
\end{aligned} 
\end{equation} 
where \(B(t)\) and \(H(t)\) denote the densities of vegetation biomass and herbivores at time \(t\), respectively, and 

\begin{equation} 
\label{eq:response} 
C(B) = \frac{B}{1+B} \left( 1-\frac{\mu B}{1+B} \right)
\end{equation} 
characterizes the defence-dependent functional response. In system~\eqref{eq:model}, vegetation grows logistically with carrying capacity \(\alpha>0\), while herbivores undergo a constant per-capita mortality rate \mbox{\(\theta>0\)} and the parameter \(m>0\) measures the conversion efficiency from consumed vegetation into herbivore biomass. The ecological interaction between the two populations is described by the nonlinear consumption function \(C(B)\). Indeed, the first factor \(B/(1+B)\) represents a classical saturating functional response accounting for feeding limitations associated with searching, handling and resource processing, whereas the second factor introduces toxin-mediated inhibition through the parameter \(\mu \in  \left[1/4,1\right]\) \cite{Li2006,Liu2008,Feng2008, Grifo2027}: for low vegetation densities, resource availability dominates the interaction and consumption increases with biomass; for sufficiently large vegetation densities, toxin accumulation reduces grazing efficiency and may eventually overcome the beneficial effects of resource availability. 

\begin{lemma}[Geometric behaviour of $C(B)$]
\label{lemma:CB}
The qualitative behaviour of the toxin-dependent consumption function \( C(B) \) depends on the toxicity parameter \(\mu\). In particular, 
\begin{enumerate} 
\item Let \( \mu \in [\frac{1}{4}, \frac{1}{2}] \), then \(C(B)\) is strictly increasing;
\item Let \( \mu \in ]\frac{1}{2}, 1] \), then \(C(B)\) possesses a unique positive critical point \( B_M=1/(2\mu-1) \), which corresponds to a global maximum, before tending to an asymptote at $C(B)=1-\mu$.
\end{enumerate} 
\end{lemma}

\begin{remark} The threshold $\mu=1/2$ separates two qualitatively different ecological regimes. For $\mu\leq1/2$, defence effects are not sufficiently strong to overcome the increase in resource availability, and the consumption rate grows monotonically with vegetation biomass. Conversely, when $\mu>1/2$, defence-induced inhibition becomes dominant at high biomass levels, generating a unimodal functional response. This transition constitutes the key mechanism responsible for the appearance of multiple coexistence states and the bifurcation scenarios investigated in the following sections. \end{remark}

\section{Equilibrium structure and local bifurcations}
\label{sec:equilibria}

The geometric properties of the functional response established in Section \ref{sec:model} provide the foundation for the equilibrium structure of system \eqref{eq:model}. In particular, the transition from a monotone to a non-monotone consumption function introduces the possibility of multiple coexistence states and distinct dynamical regimes. However, the local dynamics of defence-mediated plant-herbivore systems have been extensively investigated in previous studies \cite{Li2006, Feng2008, Liu2008, Castillo2012, Feng2013}. Therefore, our objective is not to revisit the classical qualitative analysis of the model, but rather to identify the global geometric structures that organize the dynamics. To this aim, in this section we briefly recall the main results concerning the equilibrium structure and local bifurcations that are relevant to the subsequent analysis and we refer to \cite{Liu2008,Castillo2012,Grifo2027} for more details.

\subsection{Equilibria and local stability} \label{subsec:equilibria-stability} 
Equilibria represent stationary configurations in which vegetation growth, herbivore consumption, and herbivore mortality are exactly balanced. System~\eqref{eq:model} always admits the trivial equilibrium (TE) and the herbivore-free equilibrium (HFE)

$$
[TE]: (B^\ast, H^\ast) = (0, 0)
\hspace{0.5cm}
\text{and}
\hspace{0.5cm}
[HFE]: (B^\ast, H^\ast) = (\alpha, 0)
$$
together with two possible coexistence states. These latter are determined by the intersections between the toxin-dependent consumption function \(C(B)\) and the threshold value \(\theta/m\), representing the minimum consumption level required to compensate herbivore mortality, combined with the biological admissibility requirement $0<B<\alpha$.

\begin{lemma}[Existence of $CE_-$ and $CE_+$]
\label{lemma:steady}
Besides the \(TE\) and \(HFE\) equilibria, system \eqref{eq:model} may admit two positive coexistence equilibria 

\[ 
\begin{array}{l}
CE_-:(B^\ast,H^\ast)=(B_-,H_-) = \left(\frac{m - 2 \theta - \sqrt{m(m - 4 \mu \theta)}}{2(\theta - m + \mu m)}, \frac{B_{-} \left(\alpha - B_{-} \right)}{\alpha \theta}\right),\medskip \\
CE_+: (B^\ast,H^\ast)=(B_+,H_+) = \left(\frac{m - 2 \theta + \sqrt{m(m - 4 \mu \theta)}}{2(\theta - m + \mu m)}, \frac{B_{+} \left(\alpha - B_{+} \right)}{\alpha \theta}\right).
\end{array}
\]

\noindent More precisely, the equilibrium \(CE_-\) exists and is biologically meaningful whenever $\mu<\min\{\mu_{HFE}, \mu_M\}$ or $\mu_M < \mu < \mu_{SN}$, whereas the equilibrium \(CE_+\) for $\max\{\mu_{HFE}, \mu_M\} < \mu < \mu_{SN}$, being 

\[ \mu_{HFE} = \frac{(1+\alpha)\bigl[m\alpha-\theta(1+\alpha)\bigr]} {m\alpha^{2}}, \qquad \mu_M = \frac{\alpha + 1}{2 \alpha}  \quad \text{and} \quad \mu_{SN} = \frac{m}{4\theta}. \] 
\end{lemma}

\begin{figure}[b!]
\centering
\includegraphics[width=1\textwidth]{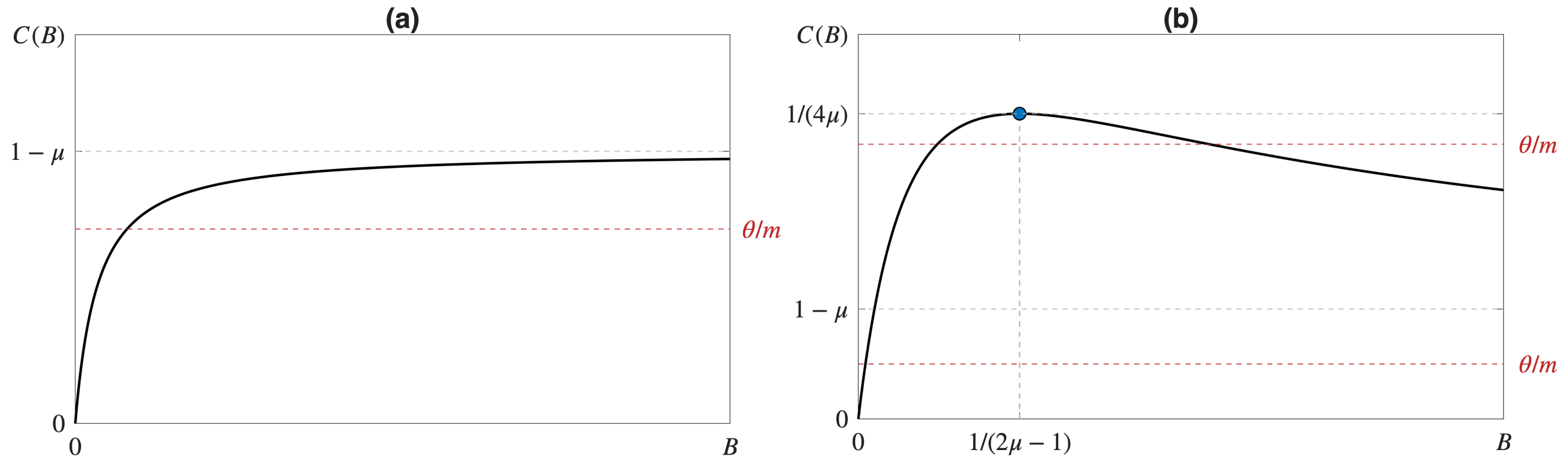}
\caption{Geometrical characterization of the toxin-dependent consumption function $C(B)$ and the nullclines of system \eqref{eq:model}. Panel (a) corresponds to $\mu\in[1/4,1/2]$, whereas panel (b) to $\mu\in(1/2,1]$.}
\label{fig:nullcline}
\end{figure}

Note that, the equation $C(B)=\theta/m$ admits at most one positive solution when \mbox{$\mu\leq 1/2$}, whereas two positive solutions may occur when $\mu>1/2$, provided that $\mu<\mu_{SN}$ \mbox{(see Figure \ref{fig:nullcline}).}

\begin{remark} 
Lemma~\ref{lemma:steady} reveals a distinction between the existence of algebraic solutions of the equation \mbox{\( mC(B)=\theta \)} and their biological admissibility. Indeed, when 

\[ 
\theta>\theta_M:=\frac{m\alpha}{2(1+\alpha)}, 
\] 
one has \( \mu_{HFE}<\mu_{SN}<\mu_M. \) Therefore, for parameter values satisfying \( \mu_{HFE}<\mu<\mu_{SN}, \) the equation \mbox{$mC(B)=\theta$} admits two distinct positive solutions $B_-$ and $B_+$, but both satisfy \( B_\pm>\alpha . \) Geometrically, this situation occurs because the horizontal level $\theta/m$ intersects the graph of $C(B)$ twice, but both intersections lie to the right of the carrying capacity threshold $B=\alpha$. As a result, the region \( \mu_{HFE}<\mu<\mu_{SN} \) contains no biologically admissible interior equilibrium despite the presence of two algebraic solutions of the coexistence equation.
\end{remark}

\noindent Ecologically, the equilibrium \(CE_-\) corresponds to a low-biomass coexistence regime, whereas \(CE_+\) is associated with a high-biomass state in which herbivores persist despite the inhibitory effects induced by plant toxicity. 

\begin{lemma}[Equilibria stability] 
\label{lem:local-classification} 
Assume that the corresponding equilibria are biologically admissible, then 

\begin{enumerate} 
\item[(i)] \(\mathrm{TE}\) is always a hyperbolic saddle,
\item[(ii)] \(\mathrm{HFE}\) is a stable node for \(\mu>\mu_{\mathrm{HFE}}\), and a hyperbolic saddle for \(\mu<\mu_{\mathrm{HFE}}\),
\item[(iii)] \(\mathrm{CE}_{-}\) is locally asymptotically stable for \(\theta>\theta_H\), and unstable for \(\theta<\theta_H\),
\item[(iv)] \(\mathrm{CE}_{+}\) is always a hyperbolic saddle, 
\end{enumerate} 
where \(\theta_H\) denotes the Hopf threshold 

\begin{equation}
\theta_H= m \frac{ \left[1+\alpha(1-\mu)\right] \left[ 1+\alpha(1+5\mu)-4\mu\right] +\left[\alpha(3\mu-1) - 1 \right] \sqrt{ \left[1+\alpha(1-\mu)\right]\left[ \alpha(1-\mu)+1+8\mu \right] } } {8\mu(1+\alpha)^2}.
\end{equation}
\end{lemma}  

As it has been already observed, the emergence of coexistence equilibria $CE_{\pm}$ is associated with the intersection between the graph of $C(B)$ and the threshold $\theta/m$, so suggesting the presence of saddle-node bifurcation mechanisms in the parameter space. 

\begin{proposition}[Saddle-node bifurcation of coexistence equilibria] 
\label{prop:SN}
Assume $\mu > 1/2$ and consider $\theta$ as bifurcation parameter. At the critical value \( \theta_{SN}= m / (4\mu), \) the coexistence equilibria $CE_-$ and $CE_+$ match at 

\begin{equation}
\label{eq:BM}
(B,H)= \left( B_M , \frac{B_M(\alpha-B_M)}{\alpha\theta_{SN}} \right). 
\end{equation}
Moreover, provided $B_M<\alpha$, system \eqref{eq:model} undergoes a codimension-one saddle-node bifurcation.
\end{proposition}

\begin{proof}
For $\mu>1/2$, Lemma~\ref{lemma:CB} implies that the consumption function $C(B)$ is unimodal and attains its unique maximum at \( B_M\) (see Figure \ref{fig:nullcline}(b)). Indeed, when \( \theta=\theta_{SN}, \) the line $\theta/m$ is tangent to $C(B)$ at $B=B_M$ and the two coexistence equilibria $CE_-$ and $CE_+$ match at \eqref{eq:BM}. For $\theta<\theta_{SN}$ the level $\theta/m$ lies below the maximum of $C(B)$ and intersects the graph twice, yielding two coexistence equilibria; for $\theta>\theta_{SN}$ no intersection exists. Hence, provided $B_M<\alpha$, the collision of $CE_-$ and $CE_+$ at $\theta=\theta_{SN}$ generates a codimension-one saddle-node bifurcation.
\end{proof}

\subsection{Hopf bifurcation and limit cycles}

The loss of stability of the coexistence equilibrium $CE_-$ is associated with the emergence of oscillatory dynamics. Let us now characterize the critical configuration at which this transition occurs.

\begin{theorem}[Hopf bifurcation at $CE_-$] 
\label{th:hopf}
Assume that $CE_-$ exists in a neighbourhood of $\theta=\theta_H$. At the critical threshold \( \theta=\theta_H \), $CE_-$ has a pair of purely imaginary eigenvalues \mbox{\( \lambda_{1,2} = \pm i \sqrt{\det J(CE_-)} \)}. Moreover, the transversality condition 

\[ \frac{d}{d\theta} \operatorname{Re}\lambda(\theta) \Big|_{\theta=\theta_H} \neq 0 \] 
is satisfied, and the first Lyapunov coefficient $l_1(\theta_H)$ is nonzero. Therefore, system~\eqref{eq:model} undergoes a nondegenerate Hopf bifurcation at $CE_-$ for $\theta=\theta_H$ with a Hopf frequency 

\[ 
\omega_0 = \sqrt{ \theta_H\,m C'(B_-(\theta_H))H_-(\theta_H)}.
\] 
\end{theorem}

\begin{proof} 
Since the equilibrium \(CE_-\) lies on the increasing branch of the toxin-dependent functional response, \( \det J(CE_-)>0 \) so that its stability is entirely determined by the sign of the Jacobian trace \( \operatorname{tr}J(CE_-) = \theta H'(B_-) \).  At the critical configuration \(\theta=\theta_H\), in which \( H'(B_-)=0 \), it follows \( \operatorname{tr}J(CE_-)=0 \). Then, the Jacobian admits a pair of purely imaginary eigenvalues \( \lambda_{1,2} = \pm i\sqrt{\det J(CE_-)} \) and the transversality condition is satisfied \cite{Grifo2027}, being 

\[
\left. \frac{d\,\mathrm{Re}\{\lambda(\theta)\}} {d\theta} \right|_{\theta=\theta_H} = \frac12 \left[ \theta H''(B_-) \frac{\partial B_-}{\partial\theta} \right]_{\theta=\theta_H} <0.
\]

Since system~\eqref{eq:model} is two-dimensional, the 2D centre manifold at the Hopf bifurcation coincides locally with the phase plane so that the normal-form computation can be performed directly in the original coordinates. To reduce the dynamics in such a manifold, let us introduce $x=B-B_{\mathrm H}$ and $v=H-H_{\mathrm H}$, being $B_{\mathrm H}:=B_-(\theta_H)$ and $H_{\mathrm H}:=H_-(\theta_H)$, and, for convenience, set 

\[ 
c_k:=C^{(k)}(B_H), \qquad k=0,1,2,3, 
\] 
with \(c_0=C(B_H)\). At $\theta=\theta_H$, the system reduces to

\[ \begin{pmatrix} 
\dot{x}\\[1mm] 
\dot{v} \end{pmatrix} = \begin{pmatrix} 0 & -\theta_H\\ m c_1 H_{\mathrm H} & 0 \end{pmatrix} 
\begin{pmatrix} 
x\\[1mm] 
v 
\end{pmatrix} + \mathcal O\bigl(\|(x,v)\|^2\bigr). 
\] 
and the corresponding Hopf frequency is given by

\[ 
\omega_0 = \sqrt{\det J(CE_-( \theta_H))} = \sqrt{ \theta_H\,m c_1 H_{\mathrm H} } >0. 
\] 
Introducing the normalized real coordinate $y = \theta_H v/ \omega_0$, the linearized system becomes 

\[ 
\dot{x}=-\omega_0y, \qquad \dot{y}=\omega_0x. 
\] 
Thus, the complex center coordinate is defined by 

\[ 
z=x+iy = B-B_{\mathrm H} + i\frac{\theta_H}{\omega_0} \left(H-H_{\mathrm H}\right)
\] 
and the vector field takes the form 

\[ 
\dot z = i\omega_0 z + F_2(z,\bar z) + F_3(z,\bar z) + O(|z|^4), 
\] 
where $F_j$ denotes the homogeneous component of degree $j$. By the Poincar\'e-Dulac normal-form theorem \cite{kuznetsov2023}, there exists a near-identity transformation 

\[ 
z=w+h_2(w,\bar w)+h_3(w,\bar w), 
\] 
where $h_j$ is homogeneous of degree $j$, which, for $\theta$ sufficiently close to $\theta_H$, eliminates all nonresonant quadratic terms. Therefore, the vector field to the centre manifold can be written as

\[ 
\dot w = (\beta(\theta)+i \Omega(\theta))w + (l_1(\theta)+i\nu_1(\theta))w|w|^2 + O(|w|^4). 
\]
with

\begin{equation}
\label{eq:betaprimo}
\beta'(\theta_H)=\frac{\theta_H}{2}
H''(B_H)B_{H}'(\theta_H) \not = 0,
\qquad
\beta(\theta_H)=0,
\qquad
\Omega(\theta_H)=\omega_0,
\end{equation}
and, at $\theta = \theta_H$, the first Lyapunov coefficient $l_1$ is given by

\begin{equation}
\label{eq:first_lyapunov}
l_1(\theta_H) = \frac{ 2(c_1^2+c_0c_2) + \alpha mH_H \left( c_0c_2^2+c_1^2c_2-c_0c_1c_3 \right) } { 4\alpha m\sqrt{H_H}\, c_0^{3/2}c_1^{3/2} }. 
\end{equation}
Note that, $l_1(\theta_H)$ may vanish only when the model parameters satisfy the degeneracy condition 

\[ 
2(c_1^2+c_0c_2) + \alpha mH_H \left( c_0c_2^2+c_1^2c_2-c_0c_1c_3 \right) =0. 
\] 
which, whenever nonempty, defines a codimension-one degeneracy locus in the parameter space. Therefore, away from this exceptional set, the equilibrium \(CE_-\) loses stability through a nondegenerate Hopf bifurcation at \mbox{\( \theta=\theta_H \)}.
\end{proof}

\begin{corollary}[Limit cycles and stability] 
\label{cor:hopf} 
Let the hypotheses of Theorem \ref{th:hopf} be satisfied. Then, $l_1(\theta_H)<0$ and every Hopf bifurcation occurring along the coexistence branch $CE_{-}$ is supercritical. Moreover, since $\beta'(\theta_H)<0$, a family of asymptotically stable periodic orbits bifurcates from $\mathrm{CE}_{-}$ for $\theta<\theta_H$. 
\end{corollary}

\begin{proof} 
At the $\theta = \theta_H$, the denominator of the Lyapunov coefficient \eqref{eq:first_lyapunov} is strictly positive so that 

\[ \operatorname{sgn}l_1(\theta_H) = \operatorname{sgn} \underbrace{\left[ 2(c_1^2+c_0c_2) + \alpha mH_H \left( c_0c_2^2+c_1^2c_2-c_0c_1c_3 \right) \right]}_{\Delta_H}. 
\] 
After some algebraic manipulations, it is possible to recast $\Delta_H$ into 

\[ 
\Delta_H = \frac{ 4B_H P(B_H,\mu) }{ (1+B_H)^7 \left[(1-\mu)B_H+\mu+1\right] }, 
\] 
where 

\[ 
P(B,\mu) = (1-\mu)^2(2\mu-1)B^4 -(\mu-1)(\mu^2+6\mu-4)B^3 +3(\mu^2+2\mu-2)B^2 +(3\mu^2-2\mu-4)B -(2\mu+1). 
\] 
For $\mu\in[1/4,1/2]$, all the coefficients of $P$ are nonpositive and its constant term is strictly negative so that $P(B,\mu) < 0$. For $\mu\in(1/2,1]$, the condition $C'(B_H)>0$ gives $0<(2\mu-1)B_H<1$ so that $P(B_H,\mu)<0$ (see the Appendix \ref{app:first_lyapunov} for more details). Therefore, $\Delta_H<0$ implies $l_1(\theta_H)<0$ and the Hopf bifurcation is supercritical. Finally, since $\beta'(\theta_H)<0$ (see \eqref{eq:betaprimo}), the stable periodic orbits occur for $\theta \lesssim \theta_H$.
\end{proof}

\noindent From an ecological viewpoint, the Hopf bifurcation highlights the transition from stable coexistence to persistent vegetation-herbivore oscillations driven by the feedback between plant recovery, herbivore grazing, and toxin-mediated inhibition.

\begin{figure}[h!] 
\centering \includegraphics[width=0.6\textwidth] {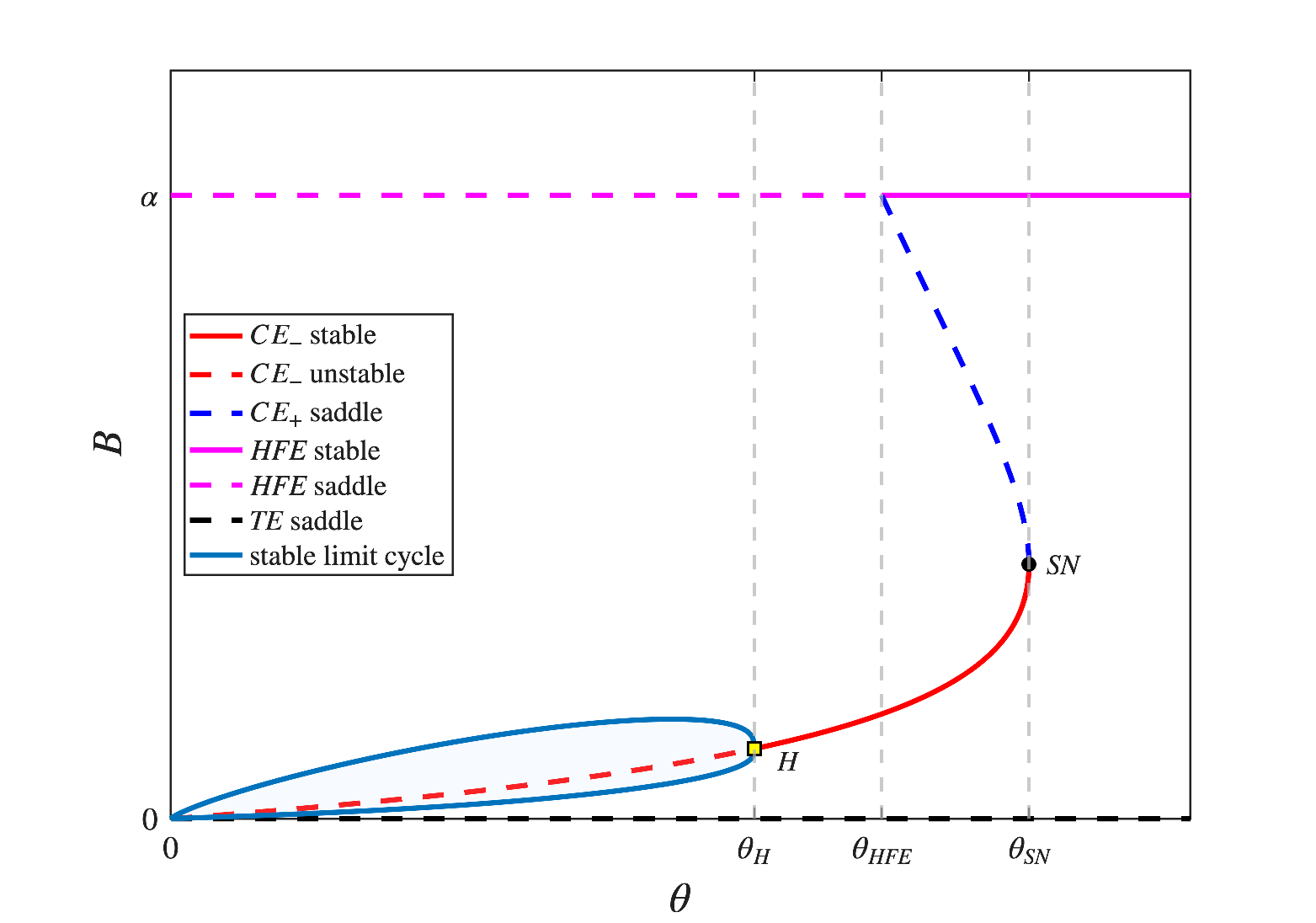} 
\caption{Qualitative bifurcation diagram with respect to the herbivore mortality parameter $\theta$. Solid and dashed curves denote stable and unstable branches, respectively. Yellow square represents the Hopf bifurcation, whereas the black circle the saddle-node one.} 
\label{fig:qualitative-bifurcation-diagram} 
\end{figure}

Before proceeding to the global phase-plane analysis, let us consider the parameter regime in which both coexistence equilibria $CE_-$ and $CE_+$ are biologically admissible. In this case, the local bifurcation structure can be summarized by the qualitative diagram shown in Figure~\ref{fig:qualitative-bifurcation-diagram}. This diagram does not represent the entire parameter space, but rather provides a schematic overview of dynamical transitions occurring as the herbivore mortality parameter $\theta$ varies. In particular, the coexistence equilibrium $CE_-$ loses stability at the Hopf threshold $\theta_H$ (yellow square), where a family of stable periodic orbits emerges (light blue solid lines), while the two coexistence branches subsequently merge at the saddle-node threshold $\theta_{SN}$ (black circle). The stability transition of the herbivore-free equilibrium $HFE$ at $\theta_{HFE}$ is also shown (magenta line). In detail, the diagram highlights four qualitatively distinct dynamical regimes: 
\begin{enumerate}
\item[i)] for $\theta<\theta_H$, the equilibrium $CE_-$ is unstable and the long-term dynamics are governed by a stable limit cycle;
\item[ii)] for $\theta_H<\theta<\theta_{HFE}$, the equilibrium $CE_-$ is asymptotically stable, providing the unique attractor of the system;
\item[iii)] for $\theta_{HFE}<\theta<\theta_{SN}$, both $CE_-$ and $HFE$ are locally asymptotically stable, leading to a bistable regime in which the asymptotic outcome depends on the initial conditions;
\item[iv)] for $\theta>\theta_{SN}$, the coexistence equilibria disappear through the saddle-node bifurcation and $HFE$ becomes the only stable equilibrium.
\end{enumerate}
Although Figure~\ref{fig:qualitative-bifurcation-diagram} illustrates only one representative bifurcation scenario, the local results obtained above completely determine the qualitative equilibrium configurations arising throughout the parameter space. Different combinations of biologically admissible equilibria and their stability properties give rise to different phase-plane organizations. We therefore proceed to a global classification in the parameter plane.

\section{Global phase-plane organization}
\label{sec:global}

The local analysis developed in Section~\ref{sec:equilibria} provides a complete characterization of the equilibria and of the bifurcation mechanisms supported by system~\eqref{eq:model}. However, local information alone is insufficient to determine the global dynamics. The long-term behaviour of the solutions is also shaped by the arrangement of invariant manifolds and by the possible occurrence of heteroclinic and homoclinic connections \cite{Xiao2015, SHAO2021}. The aim of this section is to identify the invariant structures that organize the phase portrait and to determine how they depend on the model parameters.

\subsection{Qualitative phase-plane classification} 

The existence and stability results established in Section~\ref{sec:equilibria} provide the basic ingredients required to characterize the global dynamics of system~\eqref{eq:model}. Let us investigate a qualitative classification of the phase-plane configurations, which provides the framework for the analysis of the heteroclinic and homoclinic structures discussed below.
For later convenience, let us introduce the Hopf threshold in the parameter plane by considering the inverse relation \( \mu=\mu_H(\theta), \) defined implicitly by \( \theta=\theta_H(\mu). \)

\begin{proposition}[Parameter-plane classification] 
\label{th:global}
Let system \eqref{eq:model} be defined in the biologically admissible region $\Omega= \mathbb R_+^2$. Then exactly one of the following qualitative dynamical regimes occurs:
\begin{enumerate} 
\item[(I)] If \( (\mu> \mu_{SN} \, \land \, \theta < \theta_M ) \, \lor \, ( \mu > \mu_{HFE} \, \land \, \theta> \theta_M )\), the only biologically admissible equilibria are $TE$ and $HFE$. Moreover, $TE$ is a saddle point and $HFE$ is a stable node; 

\item[(II)] If \( \max\{\mu_{H}, \mu_{HFE}\}<\mu<\mu_{SN} \, \land \, \theta<\theta_{M}\), the equilibria are \( TE,\quad HFE,\quad CE_+\) and \( CE_-\), where $TE$ and $CE_+$ are saddles whereas $HFE$ is a stable node and $CE_-$ a stable node or focus;

\item[(III)] If \( \mu_{HFE}<\mu<\mu_{H} \), the equilibria are \( TE,\quad HFE,\quad CE_+\) and \( CE_-\), where $TE$ and $CE_+$ are saddles whereas $CE_-$ is an unstable node or focus and $HFE$ is a stable node;

\item[(IV)] If \( \mu< \min\{\mu_{HFE}, \mu_H\}\), the equilibria are \( TE,\quad HFE,\quad CE_- \). Moreover,  \( TE\) and  \(HFE\) are saddles, whereas \(CE_-\) is an unstable node or focus;

\item[(V)] If \( \mu_H < \mu<\mu_{HFE}\), the equilibria are \( TE,\quad HFE,\quad CE_- \), where $TE$ and $HFE$ are saddles and $CE_-$ is a stable node or focus. 
\end{enumerate}
\end{proposition}

\begin{figure}[t!]
\centering
\includegraphics[width=0.6\textwidth]{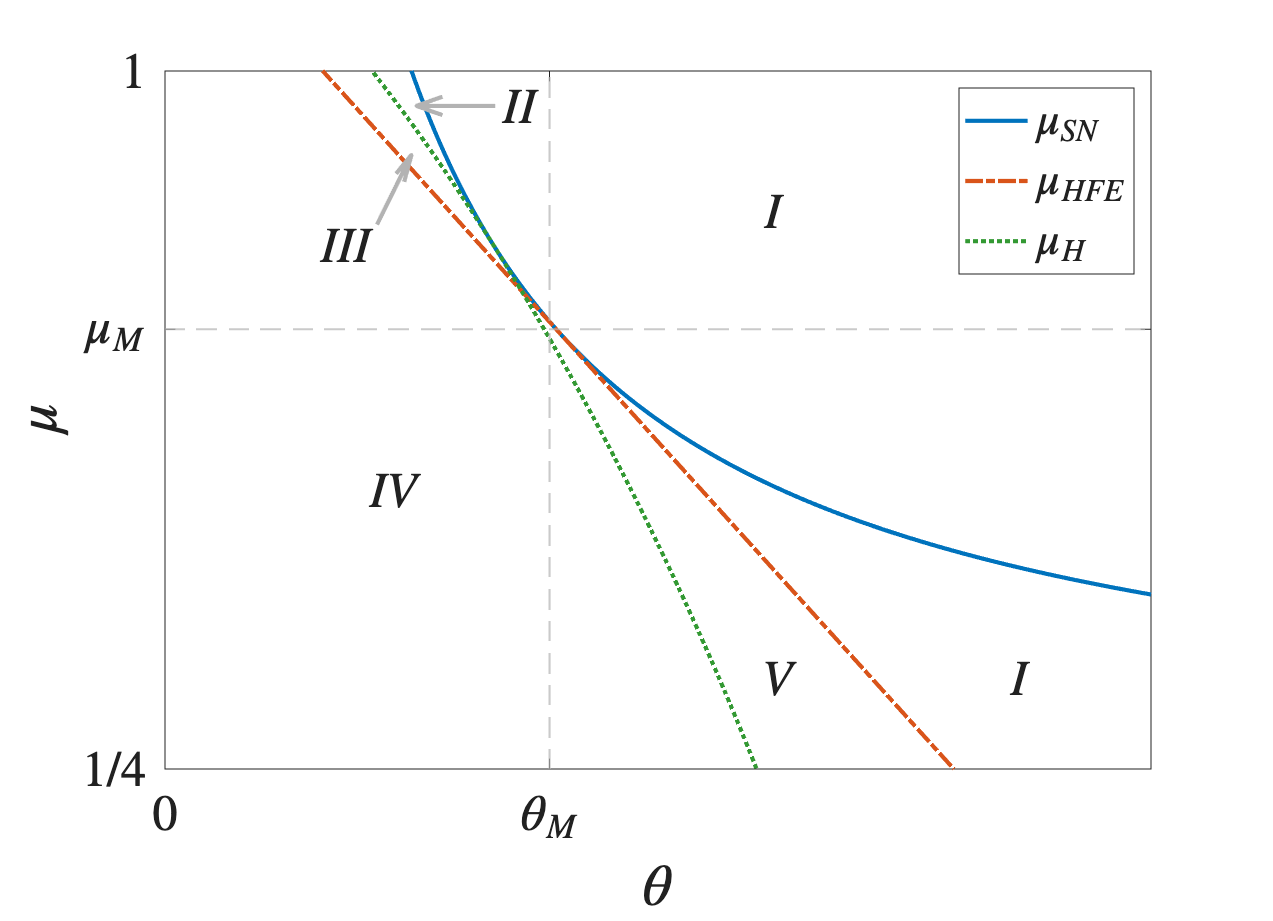}
\caption{Partition of the parameter plane \((\theta,\mu)\) induced by the thresholds \(\mu_{SN}\), \(\mu_{HFE}\), and \(\mu_H\). The different regions correspond to qualitatively distinct equilibria configurations in the phase plane of system~\eqref{eq:model}.}
\label{fig:phase_regions}
\end{figure}

\begin{proof} 
The result follows directly by combining Lemmas~\ref{lemma:steady}-\ref{lem:local-classification}. Indeed, Lemma~\ref{lemma:steady} determines the regions of the parameter space in which the coexistence equilibria \(CE_-\) and \(CE_+\) are biologically admissible, whereas Lemma~\ref{lem:local-classification} completely characterizes the stability of all equilibrium states. Consequently, the thresholds \(\mu_{HFE}\), \(\mu_{SN}\), \(\mu_H\), and \(\theta_M\) partition the parameter space into five qualitatively distinct configurations, corresponding to different arrangements of stable and unstable manifolds.
\end{proof}
The resulting decomposition is summarized in Figure~\ref{fig:phase_regions}, where each region is characterized by a specific combination of extinction and coexistence equilibria and therefore by a different set of admissible heteroclinic or homoclinic connections.

\subsubsection{Heteroclinic connections}

\begin{theorem}[Existence of heteroclinic connections] 
\label{thm:heteroclinics} 
Let us consider the system~\eqref{eq:model} in the closed positive quadrant $ \Omega =\mathbb{R}_{+}^{2} $ and assume that all the equilibria are hyperbolic. Then the following statements hold:
\begin{enumerate} 
\item[\rm (i)] The coordinate axes $\partial_B\Omega:=\{(0,H):H\geq 0\}$ and $\partial_H\Omega:=\{(B,0):B\geq 0\}$ are invariant under the flow of system~\eqref{eq:model}. Moreover, the open boundary segment $\Gamma := \{(B,0):0<B<\alpha\}$ is a heteroclinic orbit connecting the trivial equilibrium $TE$ to the herbivore-free equilibrium $HFE$. More precisely, $\Gamma \subset W^{u}(TE)\cap W^{s}(HFE)$, and 

\[ 
\lim_{t\to-\infty}\gamma(t)=TE, \qquad \lim_{t\to+\infty}\gamma(t)=HFE 
\] 
for every parametrization $\gamma$ of $\Gamma$ compatible with the flow. Hence, $TE\rightarrow HFE$ is present in every parameter region;

\item[\rm (ii)] No orbit contained in $\operatorname{int}\Omega$ can converge to $TE$ as $t\to+\infty$ so that no heteroclinic connection contained in $\operatorname{int}\Omega$ can have $TE$ as its final equilibrium. Whenever $HFE$ is a saddle, no orbit contained in $\operatorname{int}\Omega$ can converge to $HFE$ as $t\to+\infty$ so that the only heteroclinic connection terminating at $HFE$ is the boundary orbit $\Gamma$;

\item[\rm (iii)] In Region~I, the boundary connection \( TE\rightarrow HFE \) is the only heteroclinic orbit in $\Omega$;

\item[\rm (iv)] In Region~II, besides the boundary heteroclinic connection \( TE\rightarrow HFE\), there exists a heteroclinic orbit connecting the $CE_+$ to HFE. More precisely, the branch of $W^u(CE_+)$ leaving $CE_+$ in the southeastern direction is contained in the positively invariant region 
\(
\mathcal R_+:=[B_+,\alpha]\times[0,H_+] 
\) 
and converges to $HFE$ as $t\to+\infty$. Moreover, the only other heteroclinic connection between distinct equilibria that is compatible with the equilibrium configuration is \( CE_+\rightarrow CE_- \) which, in the absence of homoclinic connections generated by $CE_+$, is characterized by the branch of $W^u(CE_+)$ leaving $CE_+$ through the northwestern sector;

\item[\rm (v)] In Region III, besides the boundary heteroclinic connection \( TE\rightarrow HFE, \) there exists a heteroclinic orbit \( CE_{+}\rightarrow HFE \) as obtained in Region II. Moreover, no other heteroclinic connection between distinct equilibria contained in $\Omega$ may exist;

\item[\rm (vi)] In Region~IV, the only heteroclinic connection is \( TE\rightarrow HFE \);

\item[\rm (vii)] In Region~V, the heteroclinic connections \( TE\rightarrow HFE\) and \(HFE\rightarrow CE_{-} \) exist and are the only heteroclinic orbits between distinct equilibria contained in $\Omega$.
\end{enumerate} 
\end{theorem}

\begin{proof} 
Let 

\begin{equation} 
\label{eq:vector_field}
F(B,H) = 
\begin{pmatrix} 
B\left(1-\dfrac{B}{\alpha}\right)-mC(B)H \medskip \\
H\bigl(mC(B)-\theta\bigr) 
\end{pmatrix} 
\end{equation}
denote the vector field associated with system~\eqref{eq:model}. Since $F \in C^\infty (\Omega)$, the corresponding initial-value problem admits a unique local solution through every point of $\Omega$ \cite{Perko2001}. Therefore, let us divide the proof into three steps.\\

\noindent \textbf{Step 1: invariance of the coordinate axes and the boundary connection.} \\
The projection of system~\eqref{eq:model} to the vertical axis $B=0$ is 

\[ 
\dot B=0, \qquad \dot H=-\theta H. 
\] 
Therefore, the set \( \partial_B\Omega=\{(0,H):H\geq 0\} \) is invariant. Similarly, along the horizontal axis $H=0$ we have 

\[ 
\dot B=B\left(1-\frac{B}{\alpha}\right), \qquad \dot H=0. 
\] 
so that \( \partial_H\Omega=\{(B,0):B\geq 0\} \) is invariant too. Moreover, for $0<B<\alpha$, the scalar equation on $\partial_H\Omega$ satisfies $\dot B>0$ so that its solution with initial condition $B(0)=B_0\in(0,\alpha)$ is given by 

\[ 
B(t) = \frac{\alpha B_0e^t} {\alpha+B_0(e^t-1)}. 
\] 
Consequently, 

\[ 
\lim_{t\to-\infty}B(t)=0, \qquad \lim_{t\to+\infty}B(t)=\alpha. 
\] 
and, if \( \gamma(t):=(B(t),0), \) then 

\[ 
\lim_{t\to-\infty}\gamma(t)=TE, \qquad \lim_{t\to+\infty}\gamma(t)=HFE. 
\] 
Therefore, it follows that the boundary segment \( \Gamma \) is a heteroclinic orbit satisfying \mbox{\( \Gamma\subset W^u(TE)\cap W^s(HFE)\)} so that \( TE\rightarrow HFE \) exists for every admissible choice of the parameters.\\

\noindent \textbf{Step 2: boundary localization of the stable and unstable manifolds.}\\
At the trivial equilibrium $TE$, the eigenvalues and corresponding eigenspaces are 

\[ 
\begin{array}{l l}
\lambda_u=1, & E^u(TE)=\operatorname{span}\{(1,0)\}, \\ \medskip 
\lambda_s=-\theta, & E^s(TE)=\operatorname{span}\{(0,1)\}.
\end{array}
\] 
The invariant horizontal and vertical axes are tangent at $TE$ to $E^u(TE)$ and $E^s(TE)$, respectively. By the uniqueness of the local stable and unstable manifolds of a hyperbolic equilibrium \cite{Perko2001}, 

\[ 
W^u_{\mathrm{loc}}(TE)\cap\Omega \subseteq \{H=0\}, \qquad W^s_{\mathrm{loc}}(TE)\cap\Omega \subseteq \{B=0\}. 
\] 
We claim that no orbit contained in $\operatorname{int}\Omega$ can converge to $TE$ as $t\to+\infty$. Indeed, suppose that a solution ${\bf x}(t) = (B(t), H(t)) \in \operatorname{int}\Omega$ satisfies 

\[ 
\lim_{t\to+\infty}x(t)=TE. 
\] 
For sufficiently large $t$, the orbit must belong to $W^s_{\mathrm{loc}}(TE)$. Therefore, it must lie on the invariant axis $B=0$, contradicting \( x(t)\in\operatorname{int}\Omega. \) Hence \( W^s(TE)\cap\operatorname{int}\Omega=\varnothing \) so that no heteroclinic orbit contained in $\operatorname{int}\Omega$ can terminate at $TE$.  Note that the same argument also shows that no heteroclinic orbit from $TE$ can enter $\operatorname{int}\Omega$. Indeed, any orbit converging to $TE$ as $t\to-\infty$ must eventually belong to $W^u_{\mathrm{loc}}(TE)$, which is contained in the invariant axis $H=0$. By uniqueness of solutions, such an orbit remains on $H=0$ for all times for which it is defined. Therefore, the only heteroclinic orbit issuing from $TE$ in $\Omega$ is the boundary connection \( TE\longrightarrow HFE. \) 

Let us now consider $HFE$. Its negative eigenvalue is $\lambda_s=-1$ and the corresponding eigenspace is \mbox{\( E^s(HFE)=\operatorname{span}\{(1,0)\}. \)} Since the horizontal axis $H=0$ is invariant and tangent to $E^s(HFE)$ at $HFE$, uniqueness of the local stable manifold gives \( W^s_{\mathrm{loc}}(HFE)\cap\Omega \subseteq\{H=0\}. \) As done for $TE$, if an orbit contained in $\operatorname{int}\Omega$ converged to $HFE$ as $t\to+\infty$, then it would eventually belong to $W^s_{\mathrm{loc}}(HFE)$ and hence to $H=0$, which is impossible. Thus, whenever $HFE$ is a saddle, it follows \( W^s(HFE)\cap\operatorname{int}\Omega=\varnothing \) so that the only heteroclinic connection ending at $HFE$ is $\Gamma$.\\

\noindent \textbf{Step 3: classification in the five parameter regions.}\\ Let us now combine the previous invariant-manifold properties with the equilibrium classification established in Proposition~\ref{th:global}.\\

\noindent \emph{Region I.}\\
The only equilibria in $\Omega$ are $TE$ and $HFE$, with $TE$ a saddle and $HFE$ asymptotically stable. Since the unstable manifold of $TE$ in $\Omega$ is contained in the invariant axis $H=0$, its admissible branch generates the boundary connection $\Gamma$. There is no other saddle equilibrium from which a distinct heteroclinic orbit could originate so that \( TE\longrightarrow HFE \) is the only connection contained in $\Omega$.\\
 
\noindent \emph{Region II.}\\
In this region, $TE$ and $CE_+$ are saddles, whereas $HFE$ and $CE_-$ are asymptotically stable. As proved above, no heteroclinic orbit can terminate at $TE$ and no interior ones can exit from it (statement $(i)$). Therefore, any heteroclinic connection involving $TE$ must coincide with \( TE \rightarrow HFE \) and the only remaining saddle from which a heteroclinic orbit may issue is $CE_+$. If a branch of $W^u(CE_+)$ converges to an equilibrium, the only remaining possibilities are $HFE$ and $CE_-$. 

Let us now show that at least one branch of the unstable manifold of $CE_+$ converges to $HFE$. To this aim, let us consider the positively invariant closed rectangle \( \mathcal R_+ := [B_+,\alpha]\times[0,H_+]. \) Indeed, on the left boundary ($B=B_+$ and $0\leq H\leq H_+$) one has 

\[ 
\dot B = B_+\left(1-\frac{B_+}{\alpha}\right)-mC(B_+)H = \theta(H_+-H)\geq0,
\] 
on the right one ($B=\alpha$) it follows \( \dot B=-mC(\alpha)H\leq0 \), on the upper one ($H=H_+$ and $B_+<B\leq\alpha$) \( \dot H = H_+\bigl(mC(B)-\theta\bigr)<0\) and the lower one ($H=0$) is invariant (see Figure \ref{fig:region-II-manifolds}). 

Let $\lambda_u>0$ be the unstable eigenvalue and $v_u=(\lambda_u, mC'(B_+)H_+)^{\mathsf T}$ the associated eigenvector. Being $CE_+$ on the decreasing branch of $C$, $C'(B_+) < 0$ so that one branch of $W^u(CE_+)$ enters the interior of $\mathcal R_+$ through the southeast sector \( B>B_+\) and \(H<H_+\). Let $\gamma_{\mathrm{SE}}(t)$ denote the corresponding orbit. By the positive invariance and compactness of $\mathcal R_+$, its positive semiorbit is bounded and its $\omega$-limit set is nonempty, compact, connected, and invariant. 
Furthermore, \( \dot H <0 \) along every portion of $\gamma_{\mathrm{SE}}$ contained in $\operatorname{int}\mathcal R_+$ so that $H$ is strictly decreasing along $\gamma_{\mathrm{SE}}$, which excludes both periodic orbits and a homoclinic return to $CE_+$. Since the only equilibria contained in $\mathcal R_+$ are $CE_+$ and $HFE$, the Poincar\'e-Bendixson theorem \cite{Perko2001} yields \( \omega(\gamma_{\mathrm{SE}})=\{HFE\}. \) Hence \( \gamma_{\mathrm{SE}} \subset W^u(CE_+)\cap W^s(HFE), \) and the heteroclinic connection \( CE_+\rightarrow HFE \) exists.

Let $\gamma_{\mathrm{NW}}(t)$ denote the opposite branch of $W^u(CE_+)$, which locally leaves $CE_+$ through the northwest sector $B<B_+$ and $H>H_+$. By the boundedness of the positive solutions of system~\eqref{eq:model} \cite{Castillo2012}, its positive semiorbit is bounded. Hence $\omega(\gamma_{\mathrm{NW}})$ is nonempty, compact, connected, and invariant. We first show that the northwestern branch cannot converge to $\mathrm{HFE}$. Let \( \mathcal B(CE_-)\) and \(\mathcal B(\mathrm{HFE}) \) denote the basins of attraction of the two asymptotically stable equilibria in Region~II. These sets are nonempty, open, disjoint, and positively invariant. Moreover, since $CE_-$ is locally asymptotically stable, Lemma~\ref{lem:local-classification} gives $H'(B_-)<0$ \cite{Grifo2027}. Therefore, by Proposition~3.2 in \cite{Castillo2012}, system~\eqref{eq:model} does not admit nonconstant periodic orbit in $\operatorname{int}\Omega$. Consequently, the Poincar\'e-Bendixson theorem implies that every bounded positive orbit that does not belong to the stable manifold of $CE_+$ converges to one of the two sinks $CE_-$ and $\mathrm{HFE}$. We claim that the common boundary of the two basins in $\operatorname{int}\Omega$ is the stable manifold of $CE_+$, namely $W^s(CE_+)$, provided no homoclinic loop based at $CE_+$ exists. No point on this boundary can converge to either $CE_-$ or $\mathrm{HFE}$, since their basins are open. Its $\omega$-limit set cannot contain $TE$ (statement $(i)$) and it cannot be a periodic orbit \cite{Castillo2012}. Since $CE_+$ is the only saddle equilibrium in int$\Omega$, the Poincar\'e-Bendixson theorem gives 

\[ 
\partial\mathcal B (CE_-)\cap\operatorname{int}\Omega \subseteq W^s(CE_+), \qquad \partial\mathcal B(HFE)\cap\operatorname{int}\Omega \subseteq W^s(CE_+). 
\]
Therefore, the local stable manifold of $CE_+$ separates a sufficiently small neighbourhood of the saddle into two connected components. Since the southeastern branch of $W^u(CE_+)$ belongs to $\mathcal B (HFE)$, whereas the northwestern branch belongs to $\mathcal B (CE_-)$, every point of $W^s(CE_+)$ is an accumulation point for both basins. Therefore \(W^s(CE_+) \subseteq \partial\mathcal B(CE_-) \cap \partial\mathcal B(HFE) \) so that 

\[ 
\partial\mathcal B_{CE}\cap\operatorname{int}\Omega = \partial\mathcal B_{HFE}\cap\operatorname{int}\Omega = W^s(CE_+). 
\]
Since the southeastern branch converges to $\mathrm{HFE}$, $\gamma_{\mathrm{SE}} \subset \mathcal B(HFE)$. Therefore the northwestern branch lies on the $\mathcal B(CE_-)$ side of the stable branch and cannot converge to $\mathrm{HFE}$ so that \mbox{\( \omega(\gamma_{\mathrm{NW}}) \subseteq \{CE_-,CE_+\}. \)} If $\gamma_{\rm NW}$ returns to $CE_+$, it generates a homoclinic loop based at $CE_+$, otherwise its $\omega$-limit set cannot contain $CE_+$ and \mbox{\( \omega(\gamma_{\mathrm{NW}})=\{CE_-\}. \)} Consequently, if the northwestern branch does not form a homoclinic loop, it must converge to $CE_-$ yielding the heteroclinic connection $CE_+ \rightarrow CE_-$. However, since the existence of a homoclinic loop does not preclude the existence of such a heteroclinic connection, the contrary is not necessarily true. Note that a qualitative sketch of the above dynamics is provided in Figure \ref{fig:region-II-manifolds}.\\

\begin{figure}[t!] 
\centering 
\includegraphics[width=0.4\textwidth]{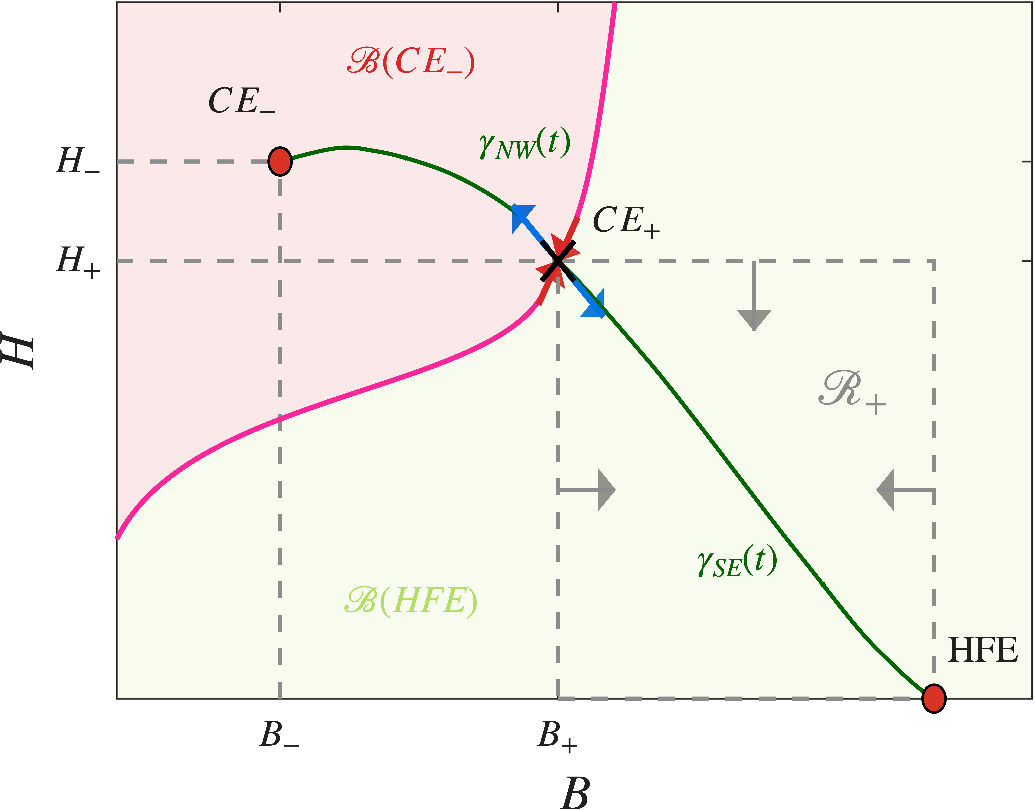} 
\caption{Geometric organization of the stable and unstable manifolds of the saddle equilibrium $\mathrm{CE}_{+}$ in Region II. The southeastern branch $\gamma_{\mathrm{SE}}$ (lower green curve) enters the positively invariant rectangle \( \mathcal{R}_{+}=[B_{+},\alpha]\times[0,H_{+}] \) via $W^{u}(\mathrm{CE}_{+})$, whereas the northwestern one $\gamma_{\mathrm{NW}}$ (upper green curve) lies on the opposite side of the stable manifold $W^{s}(\mathrm{CE}_{+})$ (magenta curve), which separates the basins of attraction of $\mathrm{CE}_{-}$ and $\mathrm{HFE}$. Blue (red) arrows indicate the direction of the eigenvectors associated with the unstable (stable) eigenspace of the linearized system at $CE_+$. Black crosses denote saddle equilibria and red circles the asymptotically stable ones.} 
\label{fig:region-II-manifolds}
\end{figure}

\noindent \emph{Region III.}\\
The argument concerning $TE$ and $CE_+$ is the same as in Region~II. In particular, the branch of $W^u(CE_+)$ entering the positively invariant region $\mathcal R_+$ converges to $HFE$, thereby generating the heteroclinic connection \( CE_+\rightarrow HFE. \) Since $CE_-$ is an unstable node or focus, its stable manifold is trivial, namely \( W^s(CE_-)=\{CE_-\} \) so that no nonconstant orbit can converge to $CE_-$ as $t\to+\infty$. Moreover, the periodic orbit $\Gamma_H$ surrounding $CE_-$ separates $CE_-$ from $HFE$ (see Theorem \ref{th:hopf}). Since trajectories of the planar system cannot cross the invariant curve $\Gamma_H$ no orbit issuing from $CE_-$ can converge to $HFE$. Thus, the heteroclinic connection \( CE_-\rightarrow HFE \) is also excluded. Thus the only heteroclinic connections contained in $\Omega$ are \( TE\rightarrow HFE \) and \(CE_+\rightarrow HFE. \)\\
  
\noindent \emph{Region IV.}\\
In this region, $TE$ and $HFE$ are saddle points, whereas $CE_-$ is a hyperbolic unstable equilibrium. The arguments of Step~2 exclude any interior heteroclinic orbit terminating at $TE$ or at $HFE$. Moreover, no interior heteroclinic orbit can issue from $TE$, since its unstable manifold in $\Omega$ is contained in $H=0$. It follows that the boundary connection \( TE\longrightarrow HFE \) is the only heteroclinic connection involving $TE$. Moreover, since $W^s(CE_-)=\{CE_-\}$, no nonconstant orbit can converge to $CE_-$ as $t\to+\infty$ so that no heteroclinic connection having $CE_-$ as terminal equilibrium can exist.
 
\noindent \emph{Region V.}\\
The equilibria $TE$ and $HFE$ are saddles, whereas $CE_-$ is locally asymptotically stable. As before, the only heteroclinic orbit issuing from $TE$ is the boundary connection \( TE\rightarrow HFE. \) 

Furthermore, since $HFE$ is a saddle whose unstable eigendirection points transversally to the invariant boundary, a biologically admissible branch of $W^u(HFE)$ may enter $\operatorname{int}\Omega$. Let $\gamma_{HFE}$ denote such a branch. Since all positive solutions of system~\eqref{eq:model} are bounded \cite{Castillo2012}, the positive semiorbit of $\gamma_{HFE}$ is bounded as well. Therefore, its $\omega$-limit set is nonempty, compact, connected, and invariant. By Step~2, no orbit contained in $\operatorname{int}\Omega$ can converge to either $TE$ or $HFE$ so that $\omega(\gamma_{HFE})$ cannot reduce to either boundary equilibrium. Since Region~V admits no periodic orbits and $CE_-$ is the only equilibrium contained in $\operatorname{int}\Omega$  \cite{Castillo2012}, the Poincar\'e-Bendixson theorem yields \( \omega(\gamma_{HFE})=\{CE_-\} \) so that \( \gamma_{HFE} \subset W^u(HFE)\cap W^s(CE_-) \) and the heteroclinic connection \mbox{\( HFE\rightarrow CE_- \)} exists.
\end{proof}

\subsubsection{Homoclinic connections}
\label{subsec:hom}

The heteroclinic classification obtained above completely determines the possible equilibrium-to-equilibrium connections supported by the system. A natural question is whether trajectories can also connect an equilibrium to itself, generating homoclinic structures. We now show that the geometric restrictions imposed by the invariant boundary and the equilibrium configuration severely limit such possibility.

\begin{theorem}[Non-existence of homoclinic orbits] 
\label{thm:no-homoclinic} 
Consider system~\eqref{eq:model} in the closed positive quadrant \( \Omega=\mathbb{R}_+^2 \) and assume that all its equilibria are hyperbolic. Then the following statements hold:

\begin{enumerate} 
\item[(i)] Neither $TE$ nor $HFE$ can be the base point of a nonconstant homoclinic orbit contained in $\Omega$;

\item[(ii)] If $CE_-$ exists, then no nonconstant homoclinic orbit to $CE_-$ is admitted;

\item[(iii)] System~\eqref{eq:model} does not admit any homoclinic orbit in Regions I, IV, and V of Proposition~\ref{th:global};

\item[(iv)] In Regions II and III, any homoclinic orbit contained in $\Omega$, if it exists, must be based at $CE_+$. Moreover, the southeastern branch of $W^u(CE_+)$ cannot belong to the homoclinic orbit. 
\end{enumerate} 
\end{theorem} 

\begin{proof} 
Since a nonconstant homoclinic orbit must converge to the same equilibrium as $t\to \pm \infty$, its base point must be a saddle. At $TE$, the stable and unstable eigenspaces are \mbox{\( E^s(TE)=\operatorname{span}\{(0,1)\}\)}, \( E^u(TE)=\operatorname{span}\{(1,0)\} \) and the coordinate axes are tangent to these eigenspaces. By uniqueness of the local stable and unstable manifolds, \( W^s_{\mathrm{loc}}(TE)\cap\Omega\subseteq\{B=0\}\) and \mbox{\( W^u_{\mathrm{loc}}(TE)\cap\Omega\subseteq\{H=0\}\)} so that a nonconstant orbit leaving $TE$ in $\Omega$ lies on the horizontal axis, whereas an orbit converging to $TE$ lies on the vertical axis. Since the two invariant axes intersect only at $TE$, they cannot form a nonconstant homoclinic loop. Thus no homoclinic orbit to $TE$ exists. 

Then, once $HFE$ is asymptotically stable, it cannot be the base point of a nonconstant homoclinic orbit. On the other hand, if it is a saddle, its local stable manifold in $\Omega$ is contained in the invariant axis $H=0$ and the corresponding unstable branch entering $\operatorname{int}\Omega$ cannot return to $HFE$, because no orbit contained in $\operatorname{int}\Omega$ can converge to $HFE$ as $t\to+\infty$ (see statement (ii) of Theorem \ref{thm:heteroclinics}). Therefore, no homoclinic orbit to $HFE$ exists. 

Whenever $CE_-$ is asymptotically stable (unstable), its unstable (stable) manifold is trivial so that \( W^s(CE_-)\cap W^u(CE_-)=\{CE_-\}, \) and hence no nonconstant homoclinic orbit to $CE_-$ can occur. 

In Regions I, IV, and V, all possible base points of a homoclinic orbit have already been excluded. Therefore, no homoclinic orbit is admitted in such regions.

Finally, in Regions II and III the only interior saddle is $CE_+$ so that any homoclinic orbit contained in $\Omega$, if present, must be based at such an equilibrium. By Theorem~\ref{thm:heteroclinics}, the branch of $W^u(CE_+)$ entering the sector $B>B_+$ and $H<H_+$ remains in $\mathcal R_+$ and converges to $HFE$. Consequently, this branch cannot return to $CE_+$ and cannot be part of a homoclinic orbit. Thus, any possible homoclinic connection to $CE_+$ would necessarily involve the opposite northwest branch of $W^u(CE_+)$. 
\end{proof}

\begin{remark} 
Assume that a homoclinic orbit $\Gamma_{\mathrm{hom}}$ based at the saddle equilibrium $CE_+$ exists. By Theorem~\ref{thm:no-homoclinic}, the orbit must leave $CE_+$ through the branch $\gamma_{NW}(t)$ via $W^u(CE_+)$. However, two geometrically distinct return configurations are possible. Indeed, the homoclinic trajectory may approach $CE_+$ either through the stable branch lying in the southwestern sector or through the stable branch lying in the northeastern one. Both configurations are compatible with the local saddle geometry and will be distinguished later when discussing the homoclinic loop. 
\end{remark}

Theorem \ref{thm:no-homoclinic} drastically restricts the possible occurrence of homoclinic connections. Having identified the unique geometric configuration compatible with such orbits, the next question is whether the corresponding stable and unstable manifolds can actually coincide for suitable parameter values. To this end, we employ a \textit{distance function} between invariant manifolds \cite{GuckenheimerHolmes1983, Perko2001, kuznetsov2023, Wiggins2003}, which provides a criterion for the existence or nonexistence of homoclinic loops associated with $CE_+$.

\begin{lemma}[Transverse section near a reference homoclinic connection] 
\label{lemma:homoclinic-section} 
Assume that, for $\eta=\eta^\ast$, system~\eqref{eq:model} admits a homoclinic orbit $\Gamma_{\rm hom}$ based at the hyperbolic saddle $CE_+(\eta^\ast)$. Suppose that $\Gamma_{\rm hom}$ leaves $CE_+(\eta^\ast)$ along the unstable branch $\gamma^u_{NW}(\eta^\ast)$ and returns to it along the selected stable branch $\gamma^s(\eta^\ast)$. Let \( q^\ast\in \Gamma_{\rm hom} \setminus\{CE_+(\eta^\ast)\}. \) Then there exist a neighbourhood $U$ of $q^\ast$, a local transverse section $\Sigma\subset U$ through $q^\ast$, and an open interval $J_0\ni\eta^\ast$ such that $\forall \eta\in J_0$: 

\begin{enumerate} 
\item[(i)] the continuation of the unstable orbit segment $\gamma^u_{NW}(\eta^\ast)$ reaching $q^\ast$ intersects $\Sigma$ transversely at a unique point $q_u(\eta)\in\Sigma\cap U$; 

\item[(ii)] the continuation of the stable orbit segment $\gamma^s(\eta^\ast)$ reaching $q^\ast$ in backward time intersects $\Sigma$ transversely at a unique point $q_s(\eta)\in\Sigma\cap U$; 

\item[(iii)] the maps \( \eta\longmapsto q_u(\eta)\) and \( \eta\longmapsto q_s(\eta) \) are continuous in $J$. If the vector field depends $C^r$-smoothly on $\eta$, then these maps are $C^r$-smooth; 

\item[(iv)] at the reference parameter, \( q_u(\eta^\ast)=q_s(\eta^\ast)=q^\ast. \)
\end{enumerate} 
\end{lemma}

\begin{figure}[h!] 
\centering 
\includegraphics[width=0.4\textwidth]{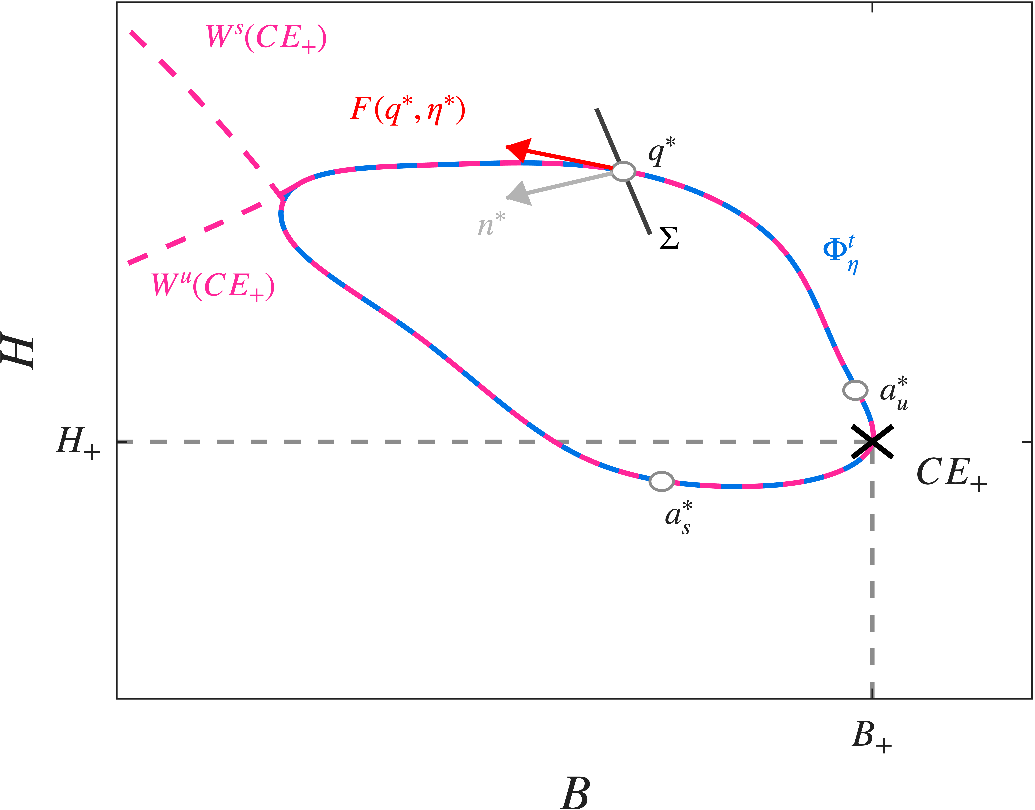} 
\caption{Local transverse section $\Sigma$ used to define the splitting function $\Delta(\eta)$ in Lemma~\ref{lemma:homoclinic-section}.} 
\label{fig:local_transverse}
\end{figure}

\begin{proof} 
Let us write system~\eqref{eq:model} as 

\[ 
\dot{x}=F(x,\eta), \qquad x=(B,H)\in\mathbb{R}^2, 
\] 
and denote by $\Phi_\eta^t$ its local flow. Then, assume that $F$ is at least $C^1$ with respect to $(x,\eta)$ and fix \( q^\ast\in \Gamma_{\rm hom}\setminus\{CE_+(\eta^\ast)\} \) so that, since the point $q^\ast$ is not an equilibrium, \( F(q^\ast,\eta^\ast)\neq 0. \) Then, let us choose a vector $n^\ast\in\mathbb{R}^2$ satisfying \( \left\langle n^\ast,F(q^\ast,\eta^\ast)\right\rangle\neq 0, \) and define 

\[ 
\Sigma := \left\{ x\in U: \left\langle n^\ast,x-q^\ast\right\rangle=0 \right\}, 
\] 
where $U$ is a sufficiently small neighbourhood of $q^\ast$. Thus, $\Sigma$ is a local section through $q^\ast$ transverse to the vector field at $\eta=\eta^\ast$ (see Figure \ref{fig:local_transverse}). By continuity of $F$, after shrinking $U$ and restricting the parameter to a sufficiently small open interval $J_0 \ni\eta^\ast$, we may assume that 

\[ 
\left\langle n^\ast,F(x,\eta)\right\rangle\neq 0  \, \forall (x,\eta)\in U\times J. 
\] 
Consequently, $\Sigma$ is transverse to the vector field for every $\eta\in J_0$. Moreover, $U$ can be chosen so that the reference homoclinic orbit intersects $U$ in a single orbit arc containing $q^\ast$. In particular, this arc intersects $\Sigma$ only at $q^\ast$. 

We first consider the unstable orbit segment. Since $\Gamma_{\rm hom}$ leaves $CE_+(\eta^\ast)$ along $\gamma^u_{NW}(\eta^\ast)$, we can choose a point 

\[ 
a_u^\ast\in W^u_{\rm loc}(CE_+(\eta^\ast)) \cap\gamma^u_{NW}(\eta^\ast) 
\] 
and a finite time $T_u^\ast>0$ such that 

\[ 
\Phi_{\eta^\ast}^{T_u^\ast}(a_u^\ast)=q^\ast. 
\]
Since $CE_+(\eta^\ast)$ is hyperbolic, hyperbolicity persists for $\eta$ sufficiently close to $\eta^\ast$. The parameter-dependent stable manifold theorem \cite{Perko2001} then provides a continuous continuation 

\[
a_u(\eta)\in W^u_{\rm loc}(CE_+(\eta)) \cap\gamma^u_{NW}(\eta), \qquad \eta\in J_0, 
\] 
such that \( a_u(\eta^\ast)=a_u^\ast. \) Define the unstable hitting function 

\[ 
G_u(t,\eta) := \left\langle n^\ast, \Phi_\eta^t(a_u(\eta))-q^\ast \right\rangle,
\] 
which, at the reference parameter, is characterized by \( G_u(T_u^\ast,\eta^\ast)=0. \) Furthermore, 

\[ 
\frac{\partial G_u}{\partial t}(T_u^\ast,\eta^\ast) = \left\langle n^\ast, F\bigl( \Phi_{\eta^\ast}^{T_u^\ast}(a_u^\ast), \eta^\ast \bigr) \right\rangle = \left\langle n^\ast,F(q^\ast,\eta^\ast) \right\rangle \neq 0. 
\] 
Therefore, by the implicit function theorem, after possibly restricting $J$, there exists a unique continuous function 

\[ 
T_u:J_0\longrightarrow\mathbb{R}, \qquad T_u(\eta^\ast)=T_u^\ast, 
\] 
such that \( G_u(T_u(\eta),\eta)=0\,\, \forall \eta\in J. \) Let us set 

\[ 
q_u(\eta) := \Phi_\eta^{T_u(\eta)}(a_u(\eta)).
\] 
It follows that \( q_u(\eta)\in\Sigma\cap U \) and \( q_u(\eta^\ast)=q^\ast. \)
The local uniqueness in the implicit function theorem, together with the choice of $U$, shows that $q_u(\eta)$ is the unique transverse intersection of the continued unstable orbit segment with $\Sigma\cap U$. 

We now consider the stable orbit segment. Since $\Gamma_{\rm hom}$ returns to $CE_+(\eta^\ast)$ along $\gamma^s(\eta^\ast)$, we can choose 

\[ 
a_s^\ast\in W^s_{\rm loc}(CE_+(\eta^\ast)) \cap\gamma^s(\eta^\ast) 
\] 
and a finite time $T_s^\ast>0$ such that 

\[ 
\Phi_{\eta^\ast}^{-T_s^\ast}(a_s^\ast)=q^\ast. 
\] 
The parameter-dependent stable manifold theorem provides a continuous continuation 

\[ 
a_s(\eta)\in W^s_{\rm loc}(CE_+(\eta)) \cap\gamma^s(\eta), \qquad \eta\in J_0, 
\] 
satisfying \( a_s(\eta^\ast)=a_s^\ast. \) Define the stable hitting function

\[ 
G_s(t,\eta) := \left\langle n^\ast, \Phi_\eta^{-t}(a_s(\eta))-q^\ast \right\rangle
\] 
so that \( G_s(T_s^\ast,\eta^\ast)=0, \) and 

\[ 
\frac{\partial G_s}{\partial t}(T_s^\ast,\eta^\ast) = - \left\langle n^\ast, F(q^\ast,\eta^\ast) \right\rangle \neq 0. 
\] 
A second application of the implicit function theorem yields a unique continuous function 

\[ 
T_s:J_0\longrightarrow\mathbb{R}, \qquad T_s(\eta^\ast)=T_s^\ast, 
\] 
such that \( G_s(T_s(\eta),\eta)=0. \) Hence, defining 

\[ 
q_s(\eta) := \Phi_\eta^{-T_s(\eta)}(a_s(\eta)), 
\] 
we obtain \( q_s(\eta)\in\Sigma\cap U\) and \( q_s(\eta^\ast)=q^\ast. \) As in the unstable case, the transversality and the choice of the neighbourhood imply that $q_s(\eta)$ is the unique intersection of the continued stable orbit segment with $\Sigma\cap U$. Note that a qualitative sketch of the above dynamics is provided in Figure~\ref{fig:local_transverse} for the case where the northwestern branch of $W^u(CE_+)$ converges to $CE_-$ through the stable branch $\gamma^s_{SW}$. Analogous phase-plane configuration arises when the convergence occurs via the stable branch $\gamma^s_{NE}$.

Therefore, we have proved that, for every $\eta\in J_0$, the selected unstable and stable orbit segments intersect $\Sigma\cap U$ transversely at the unique points $q_u(\eta)$ and $q_s(\eta)$, respectively. Moreover, 

\[ 
q_u(\eta^\ast)=q_s(\eta^\ast)=q^\ast. 
\] 
Then, the continuity of the maps 

\[ 
\eta\longmapsto q_u(\eta), \qquad \eta\longmapsto q_s(\eta) 
\] 
follows from the continuous dependence of the local invariant manifolds, the flow, and the corresponding hitting times on the parameter. Finally, if $F$ and the parameter-dependent local invariant manifolds are $C^r$ in $\eta$, the implicit function theorem in its $C^r$ version implies that $T_u$, $T_s$, $q_u$, and $q_s$ are $C^r$ as well. 
\end{proof}

\newpage

\begin{theorem}[Homoclinic bifurcation criterion] 
\label{thm:homoclinic_criterion} 
Assume that, for $\eta=\eta^\ast$, system \eqref{eq:model} admits a homoclinic orbit $\Gamma_{\rm hom}$ based at the hyperbolic saddle $CE_+(\eta^\ast)$. Suppose that $\Gamma_{\rm hom}$ leaves $CE_+(\eta^\ast)$ along the northwestern unstable branch $\gamma^u_{NW}(\eta^\ast)$ and returns to it along the selected stable branch $\gamma^s(\eta^\ast)$. Let $\Sigma$ be the transverse section and let $q_u(\eta),q_s(\eta)\in\Sigma$ be the corresponding intersection points provided by Lemma~\ref{lemma:homoclinic-section}. Let us choose a local coordinate $s:\Sigma\to\mathbb R$ and define the distance function

\[ 
\Delta(\eta) := s(q_u(\eta)) - s(q_s(\eta)). 
\] 
Then $\Delta$ is continuous in a neighbourhood of $\eta^\ast$ and 

\[ 
\Delta(\eta)=0 \iff \text{system \eqref{eq:model} admits a homoclinic orbit of the selected type}
\] 
More precisely: 
\begin{enumerate} 
\item[(i)] If $\Delta(\eta)\neq0$, then no homoclinic orbit of the selected type exists;

\item[(ii)] If there exist $\eta_1,\eta_2$ in the parameter interval under consideration such that \( \Delta(\eta_1)\Delta(\eta_2)<0, \) then there exists $\bar\eta\in(\eta_1,\eta_2)$ such that $\Delta (\bar \eta ) =0$ and system \eqref{eq:model} admits a homoclinic orbit;

\item[(iii)] If $\Delta$ is differentiable at $\eta^\ast$, \( \Delta(\eta^\ast)=0\) and \(\Delta'(\eta^\ast)\neq0, \) then the corresponding homoclinic connection is nondegenerate and defines a codimension-one homoclinic bifurcation. 
\end{enumerate} 
\end{theorem}

\begin{proof} 
By Lemma~\ref{lemma:homoclinic-section}, there exists an open interval $J_0\ni\eta^\ast$, a local transverse section $\Sigma$, and two uniquely defined transverse intersection points, namely \( q_u(\eta)\) and \(q_s(\eta), \) associated with the continuations of the selected unstable and stable orbit segments, respectively. Moreover, \mbox{\( q_u(\eta^\ast)=q_s(\eta^\ast)=q^\ast, \)} and the maps \(
\eta\mapsto q_u(\eta)\) and \( \eta\mapsto q_s(\eta) \) are continuous on $J_0$. Since the local coordinate $s:\Sigma\to\mathbb{R}$ is continuous, the splitting function \( \Delta(\eta) = s(q_u(\eta))-s(q_s(\eta)) \) is continuous on $J_0$. 

Suppose first that $\Delta(\eta)=0$, then \( s(q_u(\eta)) = s(q_s(\eta)). \) Since $s$ is a local coordinate on $\Sigma$, it is locally injective, and therefore \( q_u(\eta)=q_s(\eta)=:p. \) The point $p$ belongs simultaneously to the selected unstable and stable orbit segments so that \( p\in W^u(CE_+(\eta)) \cap W^s(CE_+(\eta)). \) By uniqueness of solutions, the two orbit segments through $p$ coincide. Consequently, the resulting orbit satisfies 

\[ 
\lim_{t\to-\infty}\gamma(t)=CE_+(\eta), \qquad \lim_{t\to+\infty}\gamma(t)=CE_+(\eta), 
\] 
and therefore defines a homoclinic orbit based at $CE_+(\eta)$. 

Conversely, assume that, for the parameter value $\eta=\eta^\ast$, system \eqref{eq:model} admits a homoclinic orbit of the selected type; namely, $\Gamma_{\mathrm{hom}}$ leaves $CE_+(\eta^\ast)$ along the northwestern unstable branch $\gamma^u_{NW}(\eta^\ast)$ and returns to it along the selected stable branch $\gamma^s(\eta^\ast)$. Since $\Gamma_{\mathrm{hom}}$ is simultaneously contained in \( W^u(CE_+(\eta^\ast))\) and \(W^s(CE_+(\eta^\ast)), \) there exists a point $q^\ast \in \Gamma_{\mathrm{hom}} \setminus \{CE_+(\eta^\ast)\}$ that belongs to both the selected unstable and stable orbit segments, namely 

\[
q^\ast \in (\Gamma_{\mathrm{hom}}(\eta^\ast) \setminus\{CE_+(\eta^\ast)\} ) \cap (W^u(CE_+(\eta^\ast)) \cap W^s(CE_+(\eta^\ast)).
\]
Applying Lemma~\ref{lemma:homoclinic-section} to the pair $(q^\ast,\eta^\ast)$, we obtain a local transverse section $\Sigma$ together with the corresponding intersection points $q_u(\eta^\ast)$ and $q_s(\eta^\ast)$. Since both orbit segments belong to the same homoclinic trajectory, they intersect $\Sigma$ at the same point. Therefore \( q_u(\eta^\ast)=q_s(\eta^\ast) \) so that 

\[ 
\Delta(\eta^\ast) = s(q_u(\eta^\ast)) - s(q_s(\eta^\ast)) = 0. 
\] 
This proves that the existence of a homoclinic orbit of the selected type implies the vanishing of the splitting function.

Having established the equivalence 

\[ 
\Delta(\eta)=0 \quad\Longleftrightarrow\quad \text{system \eqref{eq:model} admits a homoclinic orbit of the selected type}, 
\] 
the first statement follows immediately. Indeed, if $\Delta(\eta)\neq0$, then \( q_u(\eta)\neq q_s(\eta), \) and therefore the selected unstable and stable orbit segments do not coincide so that no homoclinic orbit can exist. 

Now let $\eta_1,\eta_2\in J_0$ satisfy \( \Delta(\eta_1)\Delta(\eta_2)<0. \) Since $\Delta$ is continuous on $J_0$, the Intermediate Value Theorem yields a point \( \bar\eta\in(\eta_1,\eta_2) \) such that \( \Delta(\bar\eta)=0. \) Applying again the established equivalence, it follows that system \eqref{eq:model} admits a homoclinic orbit based at $CE_+(\bar\eta)$. 

Finally, assume that \( \Delta(\eta^\ast)=0\) and \( \Delta'(\eta^\ast)\neq0. \) Then $\eta^\ast$ is a simple zero of the splitting function so that the intersection points $q_u(\eta)$ and $q_s(\eta)$ move with nonzero relative speed along the transverse section as the parameter varies. Consequently, the selected branches of $W^u(CE_+(\eta))$ and $W^s(CE_+(\eta))$ coincide only at $\eta=\eta^\ast$ and split for nearby parameter values. Therefore the corresponding homoclinic connection is nondegenerate and defines a codimension-one homoclinic bifurcation.
\end{proof}

The previous criterion characterizes the occurrence of homoclinic connections through the splitting function $\Delta$. We now establish some geometric properties that every simple homoclinic loop must satisfy.

\begin{proposition}[Geometry of simple homoclinic loops] \label{prop:geometry-homoclinic} 
Assume that system \eqref{eq:model} admits a simple homoclinic loop $\Gamma_{\rm hom}$ based at the saddle equilibrium $CE_+$. Then the following statements hold:
\begin{enumerate} 
\item[(i)] $\Gamma_{\rm hom}$ encloses $CE_-$, or equivalently $CE_- \in D_{\rm hom}$, being $D_{hom}$ the bounded connected component of $\mathbb{R}^2 \setminus \Gamma_{\rm hom}$;

\item[(ii)] Let $B_c\in(0,\alpha)$ be the unique critical point of the positive $B$-nullcline. Then $\Gamma_{\rm hom}$ intersects both regions \(\mathcal U_- = \{(B,H)\in\operatorname{int}\Omega:\ B<B_c\} \) and \( \mathcal U_+ = \{(B,H)\in\operatorname{int}\Omega:\ B>B_c\}. \) 
\end{enumerate} 
\end{proposition}

\begin{proof} 
Let $\Gamma_{\rm hom}$ be a simple homoclinic loop based at the hyperbolic saddle $CE_+$, and let $D_{\rm hom}$ denote the bounded region enclosed by $\Gamma_{\rm hom}$. We first consider the equilibria enclosed by the loop. A standard index result for planar vector fields states that the Poincar\'e index of a simple homoclinic loop based at a hyperbolic saddle, computed by a sufficiently small regularization of the loop near the saddle, is equal to one \cite{Perko2001}. Consequently, 

\[ 
\sum_{P\in D_{\rm hom}} \operatorname{ind}(P)=1, 
\] 
where the sum is taken over the equilibria $P$ lying strictly inside $ D_{\rm hom}$. Since \( \Gamma_{\rm hom}\subset\operatorname{int}\Omega, \) the bounded region $ D_{\rm hom}$ is contained in $\operatorname{int}\Omega$ so that neither $TE$ nor $HFE$ belongs to $ D_{\rm hom}$. Moreover, $CE_+$ lies on the loop and is therefore not counted among the equilibria contained strictly inside $ D_{\rm hom}$. Since $CE_-$ is a hyperbolic non-saddle equilibrium, i.e. \( \operatorname{ind}(CE_-)=1 \), and it is the only remaining equilibrium that may lie inside $ D_{\rm hom}$, the index identity then implies \( CE_-\in D_{\rm hom} \) so that the simple homoclinic loop $\Gamma_{\rm hom}$ encloses $CE_-$.

Let us now prove the second statement. Recall that the positive $B$-nullcline can be written as 

\[ 
H=N(B), \qquad N(B):= \frac{(1-B/\alpha)(1+B)^2} {m[1+(1-\mu)B]}. 
\] 
Consider the Dulac function \cite{Perko2001} 

\[ 
D(B,H):=\frac{1}{C(B)H}, \qquad (B,H)\in\operatorname{int}\Omega. 
\] 
Since $C(B)>0$ for $B>0$, the function $D$ is continuously differentiable on $\operatorname{int}\Omega$. Consequently, 

\[ 
\nabla\cdot\bigl(D(B,H)F(B,H)\bigr) = \frac{\partial}{\partial B} \left( \frac{mN(B)}{H}-m \right) + \frac{\partial}{\partial H} \left( m-\frac{\theta}{C(B)} \right) = \frac{mN'(B)}{H},
\] 
being $F(B, H)$ the vector field \eqref{eq:vector_field}.
A direct differentiation gives 

\begin{equation}
\label{eq:NB}
N'(B) = \frac{1+B} {m\alpha[1+(1-\mu)B]^2} Q(B), 
\end{equation}
where 

\[ Q(B) = \alpha(1+\mu) + \bigl[\alpha(1-\mu)-2(1+\mu)\bigr]B - 3(1-\mu)B^2. 
\] 
Since the prefactor in \eqref{eq:NB} is strictly positive, $\operatorname{sgn} N'(B) = \operatorname{sgn} Q(B)$. For $\mu<1$, the constant term of $Q$ is positive and the product of its roots is negative so that it possesses exactly one positive root. Since \( Q(0)=\alpha(1+\mu)>0\) and \(Q(\alpha) = -\alpha(1+\alpha)\bigl[1+\alpha(1-\mu)\bigr] <0, \) the unique positive root 

\[ 
B_c = \frac{ \alpha(1-\mu)-2(1+\mu) + \sqrt{ [\alpha(1-\mu)-2(1+\mu)]^2 + 12\alpha(1-\mu)(1+\mu) } } {6(1-\mu)} 
\]
belongs to $(0,\alpha)$. Moreover, $Q$ is a strictly concave quadratic polynomial, i.e.

\[ 
\begin{array}{c c c}
N'(B)>0 & \text{for} & 0<B<B_c, \medskip \\
N'(B)<0 & \text{for} & B_c<B<\alpha. 
\end{array}
\] 
Suppose, by contradiction, that $\Gamma_{\rm hom}$ were entirely contained in the region $\mathcal U_-$, region in which \( \nabla\cdot\bigl(DF\bigr) >0. \) Since $\mathcal U_-$ is simply connected, the Bendixson-Dulac criterion excludes a closed invariant curve formed by regular orbit segments, including a homoclinic loop, entirely contained in $\mathcal U_-$ \cite{Perko2001}. Similarly, if $\Gamma_{\rm hom}$ were entirely contained in $\mathcal U_+$, then \mbox{$\nabla\cdot\bigl(DF\bigr) <0 $} throughout $\mathcal U_+$. A second application of the Bendixson--Dulac criterion gives a contradiction. Therefore, $\Gamma_{\rm hom}$ cannot be contained entirely on either side of the line $B=B_c$ so that it intersects both regions $\mathcal U_\pm$.
\end{proof}

The geometric restrictions established above provide only a partial description of the homoclinic configuration. To understand the dynamical role played by a simple homoclinic loop, it is necessary to investigate the invariant regions generated by the loop and their relation with the other invariant sets of the system. The next results show how these geometric configurations determine the basin structure and the long-term dynamics in Regions II and III.

\begin{lemma}[Relative position of $B_+$ and $B_c$] 
\label{lem:Bplus-Bc} 
Assume that the coexistence equilibrium $CE_+$ exists and let $B_c\in(0,\alpha)$ denote the unique critical point of the positive $B$-nullcline. Then the following statements hold. 

\begin{enumerate} 
\item[(i)] If the parameters belong to Region~II, then \( B_+ >  B_c. \) 

\item[(ii)] If the parameters belong to Region~III, then 

\[ 
B_+>B_c 
\quad \iff \quad
\alpha< \alpha_c := \frac{B_+\left[2(1+\mu)+3(1-\mu)B_+\right]} {1+\mu+(1-\mu)B_+}. \]
\end{enumerate} 
\end{lemma}

\begin{proof} 
Let us recall from Proposition~\ref{prop:geometry-homoclinic} that \( \operatorname{sgn}N'(B) = \operatorname{sgn}Q(B). \) Moreover, $Q$ admits a unique positive root $B_c$ so that \( Q(B)>0\) for $B<B_c$ and \(Q(B)<0\) for $B>B_c$. 

Let us first consider Region~II. By Lemma~\ref{lem:local-classification}, the equilibrium $CE_-$ is asymptotically stable. Since \( \operatorname{tr}J(CE_-) = \theta N'(B_-) \), stability implies \( N'(B_-)<0. \) Hence \( B_->B_c \) and, since $B_-<B_+$, it follows that \( B_c<B_-<B_+. \) 

Let now the parameters belong to Region~III. Since $B_c$ is the unique positive root of $Q$, one has \( B_+>B_c\) iff \(Q(B_+)<0. \) Evaluating $Q$ at $B_+$ yields

\[ 
Q(B_+) = \alpha\bigl[1+\mu+(1-\mu)B_+\bigr] - B_+\bigl[2(1+\mu)+3(1-\mu)B_+\bigr],
\] 
so that 

\[ 
Q(B_+)<0 
\iff
\alpha< \frac{B_+\bigl[2(1+\mu)+3(1-\mu)B_+\bigr]} {1+\mu+(1-\mu)B_+}. 
\] 
Combining the previous equivalences, we conclude that 

\[ 
B_+>B_c \iff \alpha< \frac{B_+\bigl[2(1+\mu)+3(1-\mu)B_+\bigr]} {1+\mu+(1-\mu)B_+}. 
\] 
\end{proof}

\subsubsection*{Region II}

\begin{theorem}[Nonexistence of homoclinic loops in Region II] 
\label{prop:no-homoclinic-region-II} 
Assume that the parameters belong to Region~II. Then system~\eqref{eq:model} does not admit a simple homoclinic loop based at the saddle equilibrium $CE_+$. 
\end{theorem} 

\begin{figure}[h!] 
\centering 
\includegraphics[width=0.4\textwidth]{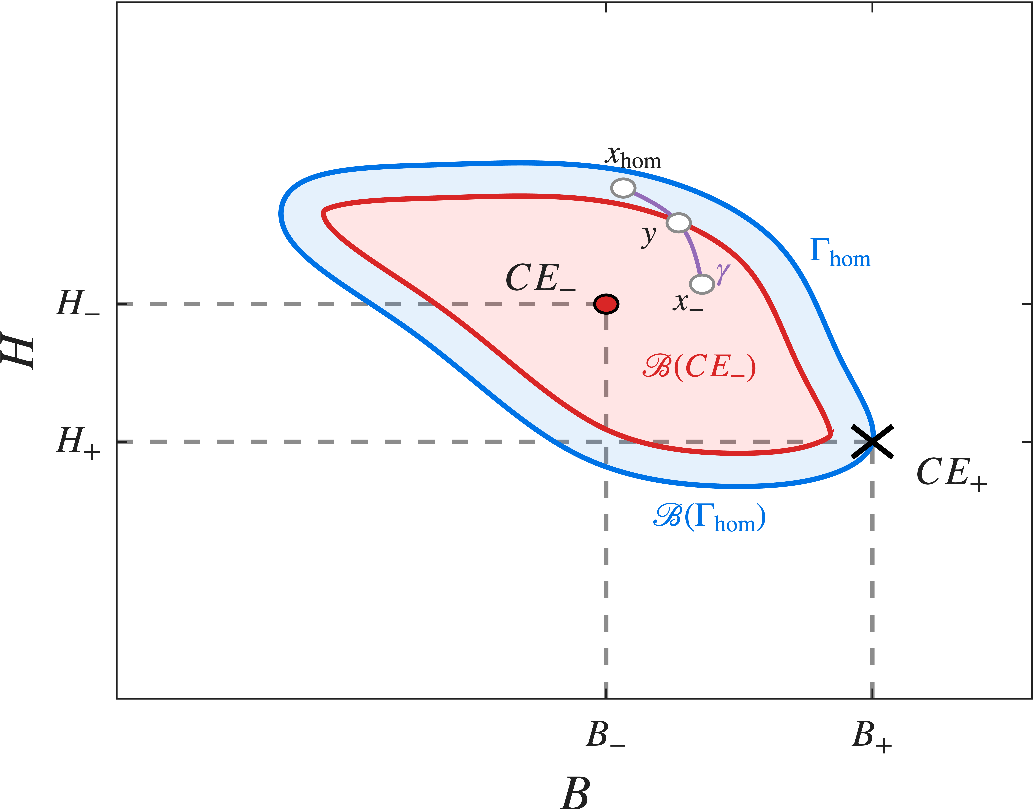} 
\caption{Geometric configuration associated with a hypothetical simple homoclinic loop $\Gamma_{\mathrm{hom}}$, supposed in Theorem~\ref{prop:no-homoclinic-region-II}, based at the saddle equilibrium $CE_{+}$ in Region~II.} 
\label{fig:hom_regionII}
\end{figure}

\begin{proof} 
Assume, by contradiction, that system~\eqref{eq:model} admits a simple homoclinic loop $\Gamma_{\mathrm{hom}}$ based at $CE_+$, and let $D_{\mathrm{hom}}$ denote the bounded connected component of $\mathbb{R}^2\setminus\Gamma_{\mathrm{hom}}$. By Proposition~\ref{prop:geometry-homoclinic}, the loop encloses the coexistence stable equilibrium $CE_-$, so that \( CE_-\in D_{\mathrm{hom}}. \) Let $\lambda_s<0<\lambda_u$ denote the eigenvalues of the Jacobian matrix at $CE_+$. The saddle quantity \cite{Perko2001} associated with $CE_+$ is 

\begin{equation}
\label{eq:sigma+}
\sigma_+:=\lambda_s+\lambda_u =\operatorname{tr}J(CE_+) = \theta N'(B_+). 
\end{equation}
By Lemma~\ref{lem:Bplus-Bc}, one has \( B_+>B_c \) throughout Region~II. Since $B_c$ is the unique critical point of $N$ and \( N'(B)<0\) for $B>B_c$, from \eqref{eq:sigma+} it follows that \( \sigma_+ <0. \) The classical local theory of planar homoclinic loops then implies that $\Gamma_{\mathrm{hom}}$ is attracting from the side on which the local Poincar\'e return map is defined \cite{GuckenheimerHolmes1983, Perko2001}. Since $\Gamma_{\mathrm{hom}}$ is a simple loop enclosing $D_{\mathrm{hom}}$, this return map is defined on a one-sided transverse section contained in $D_{\mathrm{hom}}$. Consequently, there exists a neighbourhood $\mathcal{U}_{\mathrm{hom}}$ of $\Gamma_{\mathrm{hom}}$ in $\overline{D_{\mathrm{hom}}}$ such that 

\[ 
\operatorname{dist} \bigl(\Phi_t(x),\Gamma_{\mathrm{hom}}\bigr) \longrightarrow 0 \qquad\text{as }t\to+\infty 
\] 
for every \( x\in\mathcal{U}_{\mathrm{hom}}\cap D_{\mathrm{hom}}. \) On the other hand, since $CE_-$ is locally asymptotically stable, there exists a positively invariant neighbourhood $\mathcal{U}_-$ of $CE_-$ such that 

\[ 
\operatorname{dist} \bigl(\Phi_t(x), CE_- \bigr) \longrightarrow 0 \qquad\text{as }t\to+\infty 
\] 
for every $x\in\mathcal{U}_-$. Let us choose \( x_-\in\mathcal{U}_-\) and \(x_{\mathrm{hom}}\in \mathcal{U}_{\mathrm{hom}}\cap D_{\mathrm{hom}}. \) Since $D_{\mathrm{hom}}$ is an open connected subset of $\mathbb{R}^2$, it is path connected. Hence there exists a continuous path \( \gamma\colon[0,1]\longrightarrow D_{\mathrm{hom}} \) such that \( \gamma(0)=x_-\) and \( \gamma(1)=x_{\mathrm{hom}}\) (see Figure \ref{fig:hom_regionII}). Let 

\[ 
\mathcal{B}(CE) := \left\{ x\in D_{\mathrm{hom}}: \Phi_t(x)\longrightarrow CE_- \text{ as }t\to+\infty \right\}. 
\] 
Since $CE_-$ is asymptotically stable, $\mathcal{B}(CE)$ is open such that $\mathcal{B}(CE) \subseteq D_{\mathrm{hom}}$. Its boundary \( \partial \mathcal{B}(CE) \) is nonempty, because $\gamma(0)\in\mathcal{B}(CE)$ whereas $\gamma(1)\notin\mathcal{B}(CE)$. Let us choose \( y\in\gamma([0,1]) \cap\partial \mathcal{B}(CE). \) The basin boundary is positively invariant. Since $D_{\mathrm{hom}}$ is bounded and positively invariant, the positive semiorbit of $y$ is bounded, and therefore $\omega(y)$ is nonempty, compact, connected, and invariant. The set $\omega(y)$ cannot contain $CE_-$. Indeed, if $CE_-\in\omega(y)$, then the local asymptotic stability of $CE_-$ would imply \( \Phi_t(y)\rightarrow CE_-, \) and hence $y\in\mathcal{B}(CE)$, contradicting \( y\in\partial \mathcal{B}(CE). \) Moreover, $\omega(y)$ cannot coincide with $\Gamma_{\mathrm{hom}}$. Indeed, the basin of the attracting homoclinic loop is open on the $D_{\mathrm{hom}}$ side, whereas $y$ belongs to the invariant boundary separating that basin from $\mathcal{B}(CE)$. Thus $\omega(y)$ contains no equilibrium and, by the Poincar\'e-Bendixson theorem, $\omega(y)$ is a nonconstant periodic orbit contained in $D_{\mathrm{hom}}$. This contradicts the nonexistence of nonconstant periodic orbits in Region~II established in Proposition~3.2 of~\cite{Castillo2012}. Therefore, system~\eqref{eq:model} cannot admit a simple homoclinic loop based at $CE_+$ in Region~II. A qualitative sketch of the above structures is provided in Figure~\ref{fig:hom_regionII}.
\end{proof}

\subsubsection*{Region III}

\begin{lemma}[Relative position of $\Gamma_H$ and $\Gamma_{\mathrm{hom}}$] 
\label{lem:relative_position} 
Assume that system~\eqref{eq:model} admits a simple homoclinic loop $\Gamma_{\mathrm{hom}}$ based at $CE_+$ in Region~III. Let $D_{\mathrm{hom}}$ and $D_H$ denote the bounded domains enclosed by $\Gamma_{\mathrm{hom}}$ and $\Gamma_H$, respectively. Then exactly one of the following alternatives holds: 

\begin{enumerate}
\item[(i)] If $CE_+\notin D_H$, then \( \overline{D_H}\subset D_{\mathrm{hom}}, \) that is, the stable periodic orbit $\Gamma_H$ is enclosed by the homoclinic loop;
\item[(ii)] If $CE_+\in D_H$, then \( \overline{D_{\mathrm{hom}}}\subset D_H, \) that is, the homoclinic loop is enclosed by the stable periodic orbit $\Gamma_H$. 
\end{enumerate} 
\end{lemma}

\begin{proof} 
By Proposition~\ref{prop:geometry-homoclinic}, \( CE_-\in D_{\mathrm{hom}}. \) Moreover, since $\Gamma_H$ surrounds $CE_-$, one also has \mbox{\( CE_-\in D_H. \)} Since both $\Gamma_H$ and $\Gamma_{\mathrm{hom}}\setminus\{CE_+\}$ are trajectories of system~\eqref{eq:model}, uniqueness of solutions implies that they cannot intersect, namely \( \Gamma_H\cap\Gamma_{\mathrm{hom}}=\varnothing . \) Consequently, the corresponding bounded domains must be nested. If $CE_+\notin D_H$, then the homoclinic loop cannot be contained in $D_H$ because $CE_+\in\Gamma_{\mathrm{hom}}$ so that \( \overline{D_H}\subset D_{\mathrm{hom}}. \) Conversely, if $CE_+\in D_H$, then the homoclinic loop is entirely contained in $D_H$, yielding \( \overline{D_{\mathrm{hom}}}\subset D_H. \) The two alternatives are mutually exclusive. 
\end{proof}

\begin{theorem}[Dynamical role of simple homoclinic loops in Region III] 
\label{th:homoclinic_regionIII}
Assume that the parameters belong to Region III and that system~\eqref{eq:model} admits a simple homoclinic loop $\Gamma_{\rm hom}$ based at the saddle equilibrium $CE_+$. Let $\Gamma_H$ denote the asymptotically stable periodic orbit generated through the supercritical Hopf bifurcation of $CE_-$. Let $D_{\rm hom}$ and $D_H$ denote the bounded domains enclosed by $\Gamma_{\rm hom}$ and $\Gamma_H$, respectively. Then the following statements hold:

\begin{enumerate} 
\item[(i)] Assume that \( D_H\subset D_{\rm hom}. \) If \( 
\alpha< \alpha_c, \) then $\Gamma_{\rm hom}$ is asymptotically stable from the $D_{\rm hom}$ side. Moreover, $\Gamma_H$ and $\Gamma_{\rm hom}$ are distinct attractors contained in the same positively invariant domain $D_{\rm hom}$ so that there exists a nonempty compact invariant set \( K\subset D_{\rm hom} \) such that 

\[ 
K\cap\Gamma_H=\varnothing, \qquad K\cap\Gamma_{\rm hom}=\varnothing, 
\] 
and \( K \subset \partial\mathcal B(\Gamma_H) \cap \partial\mathcal B(\Gamma_{\rm hom}). \) In particular, system~\eqref{eq:model} possesses an additional invariant structure inside $D_{\rm hom}$;

\item[(ii)] Assume that \( D_H\subset D_{\rm hom}. \) If \( \alpha> \alpha_c, \) then $\Gamma_{\rm hom}$ is repelling from the $D_{\rm hom}$ side. Furthermore, if $\Gamma_H$ is the only periodic orbit contained in $D_{\rm hom}$, then \( D_{\rm hom}\setminus\{CE_-\} \subset \mathcal B(\Gamma_H), \) that is, 

\[ 
\lim_{t\to+\infty} \operatorname{dist} \bigl( \Phi_t(x), \Gamma_H \bigr) =0 \qquad \forall x\in D_{\rm hom}\setminus\{CE_-\};
\] 

\item[(iii)] Assume that \( D_{\rm hom}\subset D_H \), \( \alpha < \alpha_c \) and that $\Gamma_{\rm hom}$ is the only periodic invariant set contained in $D_{\rm hom}$. Then, $\Gamma_{\rm hom}$ is asymptotically stable from the $D_{\rm hom}$ side and \( D_{\rm hom}\setminus\{CE_-\} \subset \mathcal B(\Gamma_{\rm hom}), \) that is, 

\[ 
\lim_{t\to+\infty} {\rm dist} \bigl( \Phi_t(x), \Gamma_{\rm hom} \bigr)=0 \qquad \forall\, x\in D_{\rm hom}\setminus\{CE_-\}. 
\] 

Moreover, if $\Gamma_H$ is the only periodic orbit contained in $D_H\setminus\overline{D_{\rm hom}}$, then \mbox{\( D_H\setminus \overline{D_{\rm hom}} \subset \mathcal B(\Gamma_H) \cup W^s(CE_+), \)} that is, 

\[ 
\lim_{t\to+\infty} {\rm dist} \bigl( \Phi_t(x), \Gamma_H \bigr)=0 \qquad \forall x\in D_H\setminus \{\overline{D_{\rm hom}} \cup W^s(CE_+)\}. 
\]
Finally, the coexistence of the two attractors $\Gamma_H$ and $HFE$ implies the existence of a nonempty compact invariant set \( K\subset \mathbb R_+^2 \setminus D_H \) contained in the common boundary of the corresponding attraction basins, \( K \subset \partial\mathcal B(\Gamma_H) \cap \partial\mathcal B(HFE). \) In particular, system~\eqref{eq:model} possesses an additional invariant structure separating the asymptotic regimes associated with $\Gamma_H$ and $HFE$;

\item[(iv)] Assume that \( D_{\rm hom}\subset D_H. \) If \( \alpha > \alpha_c, \) then $\Gamma_{\rm hom}$ is repelling from the $D_{\rm hom}$ side. Moreover, system~\eqref{eq:model} possesses at least one compact invariant set $K \subset D_{\rm hom} \setminus \{CE_-\}$. In particular, $D_{\rm hom}$ contains at least one periodic orbit $K$ distinct from $\Gamma_H$. 
\end{enumerate} 
The degenerate case \( \alpha = \alpha_C \) cannot be determined from the saddle quantity alone and requires a higher-order analysis of the associated Poincar\'e return map. \end{theorem}

\begin{proof}
By Lemma~\ref{lem:relative_position}, the periodic orbit $\Gamma_H$ and the homoclinic loop $\Gamma_{\rm hom}$ are disjoint, and exactly one of the following two alternatives holds: \( D_H\subset D_{\rm hom}\) or \(D_{\rm hom}\subset D_H\). Moreover, by Proposition~\eqref{prop:geometry-homoclinic}, \( CE_-\in D_{\rm hom}. \) Since $\Gamma_{\rm hom}$ is a simple closed invariant curve, uniqueness of solutions implies that no trajectory can cross it so that $D_{\rm hom}$ is positively invariant. Let $\lambda_s<0<\lambda_u$ be the eigenvalues of the Jacobian matrix at $CE_+$ and let \( \sigma_+:=\lambda_s+\lambda_u \) denote the associated saddle quantity. By Lemma~\ref{lem:Bplus-Bc}, it follows that 

\[
\sigma_+\gtrless0 \quad\Longleftrightarrow\quad \alpha\gtrless\alpha_c.
\] 
Consequently, by the classical local theory of planar homoclinic loops \cite{Perko2001}, $\alpha < \alpha_c$ implies that $\Gamma_{\rm hom}$ is asymptotically stable from the $D_{\rm hom}$ side. Conversely, when $\alpha > \alpha_c$, $\Gamma_{\rm hom}$ is repelling from the $D_{\rm hom}$ side. Let us now consider separately the four configurations appearing in the statement. \medskip 

\noindent \textit{Proof of statement~{\rm (i)}.} \\
Assume that \( D_H\subset D_{\rm hom}\) and \( \alpha<\alpha_c \) so that $\Gamma_{\rm hom}$ is asymptotically stable from the $D_{\rm hom}$ side and there exists a nonempty open set \( U_H\subset D_{\rm hom} \) such that 

\[ 
\lim_{t\to+\infty} \operatorname{dist}\bigl(\Phi_t(x),\Gamma_H\bigr)=0 \qquad \text{for every }x\in U_H. 
\] 
Similarly, the one-sided asymptotic stability of $\Gamma_{\rm hom}$ provides a nonempty open set \( U_{\rm hom}\subset D_{\rm hom} \) such that 

\[ 
\lim_{t\to+\infty} \operatorname{dist}\bigl(\Phi_t(x),\Gamma_{\rm hom}\bigr)=0 \qquad \text{for every }x\in U_{\rm hom}. 
\] 
Thus the relative basins \( \mathcal B (\Gamma_H)\) and \( \mathcal B_{D_{\rm hom}}(\Gamma_{\rm hom}) := \mathcal B (\Gamma_{\rm hom}) \cap D_{\rm hom}\) are nonempty, open, and disjoint subsets of the connected domain $D_{\rm hom}$. Consequently, their relative boundaries in $D_{\rm hom}$ are nonempty so that let us define \( S:= \partial \mathcal B_{D_{\rm hom}}(\Gamma_H) \cap \partial \mathcal B (\Gamma_{\rm hom}). \) Then $S$ is nonempty, closed relative to $D_{\rm hom}$, and invariant under the flow. Indeed, let us choose $y\in S$. Since $\overline{D_{\rm hom}}$ is compact and positively invariant, the positive semiorbit of $y$ is bounded and \( K:=\omega(y) \) is a nonempty compact invariant set contained in $\overline{D_{\rm hom}}$ and it cannot intersect $\Gamma_H$. Indeed, if $\omega(y)\cap\Gamma_H\neq\varnothing$, the asymptotic stability of $\Gamma_H$ would imply that \( \omega(y)=\Gamma_H \) and hence that $y$ belongs to the basin of $\Gamma_H$, contradicting the fact that $y$ belongs to its boundary. The same argument applies to $\Gamma_{\rm hom}$, using its one-sided asymptotic stability. Therefore, 

\[ 
K\cap\Gamma_H=\varnothing, \qquad K\cap\Gamma_{\rm hom}=\varnothing , \qquad \text{with} \qquad K\subset \partial\mathcal B(\Gamma_H) \cap \partial\mathcal B(\Gamma_{\rm hom}). 
\] 
Thus the coexistence of the two attractors forces the existence of an additional compact invariant structure in $D_{\rm hom}$. \\

\noindent \textit{Proof of statement~{\rm (ii)}.} \\
Assume that \( D_H\subset D_{\rm hom}\) and \(\alpha>\alpha_c \) so that $\Gamma_{\rm hom}$ is repelling from the $D_{\rm hom}$ side. Let \mbox{\( x\in D_{\rm hom}\setminus\{CE_-\}. \)} Since $D_{\rm hom}$ is bounded and positively invariant, the positive semiorbit of $x$ is bounded. It follows that $\omega(x)$ is nonempty, compact, connected, and invariant. In Region III, the equilibrium $CE_-$ is a hyperbolic repeller so that \( W^s(CE_-)=\{CE_-\}, \) and hence \( \omega(x)\neq\{CE_-\}. \) Furthermore, the local stable branches of the saddle $CE_+$ forming the homoclinic configuration belong to the boundary of $D_{\rm hom}$ so that \( W^s(CE_+)\cap D_{\rm hom}=\varnothing, \) and no point of $D_{\rm hom}$ can converge to $CE_+$. Finally, since $\Gamma_{\rm hom}$ is repelling from the $D_{\rm hom}$ side, it cannot be the $\omega$-limit set of a point in $D_{\rm hom}$. Therefore $\omega(x)$ contains no equilibrium and cannot coincide with $\Gamma_{\rm hom}$. By the Poincar\'e-Bendixson theorem, $\omega(x)$ is a periodic orbit contained in $D_{\rm hom}$. If $\Gamma_H$ is the only periodic orbit contained in $D_{\rm hom}$, then necessarily \( \omega(x)=\Gamma_H. \) Hence 

\[ 
\lim_{t\to+\infty} \operatorname{dist}\bigl(\Phi_t(x),\Gamma_H\bigr)=0 \qquad \text{for every } x\in D_{\rm hom}\setminus\{CE_-\}. 
\] \\

\noindent \textit{Proof of statement~{\rm (iii)}.}\\
Assume that \( D_{\rm hom}\subset D_H\) and \(\alpha<\alpha_c \) so that $\Gamma_{\rm hom}$ is asymptotically stable from the $D_{\rm hom}$ side. Assume that $\Gamma_{\rm hom}$ is the only periodic invariant set contained in $\overline{D_{\rm hom}}$ and let \( x\in D_{\rm hom}\setminus\{CE_-\}. \) As before, $D_{\rm hom}$ is bounded and positively invariant, so that $\omega(x)$ is nonempty, compact, connected, and invariant. Since $CE_-$ is a hyperbolic repeller, \( \omega(x)\neq\{CE_-\}, \) and, since \( W^s(CE_+)\cap D_{\rm hom}=\varnothing, \) the set $\omega(x)$ cannot contain $CE_+$. If $\omega(x)$ contains no equilibrium, the Poincar\'e-Bendixson theorem implies that it is either a periodic orbit or a finite union of equilibria and connecting trajectories. The latter alternative is excluded because the only saddle on $\partial D_{\rm hom}$ is $CE_+$ and the only corresponding invariant graphic is the homoclinic loop itself. By the uniqueness assumption on periodic invariant sets, one obtains \( \omega(x)=\Gamma_{\rm hom} \) so that 

\[ 
\lim_{t\to+\infty} \operatorname{dist} \bigl(\Phi_t(x),\Gamma_{\rm hom}\bigr)=0 \qquad \text{for every } x\in D_{\rm hom}\setminus\{CE_-\}. 
\] 
Consider now the open annular region \( A:= D_H\setminus\overline{D_{\rm hom}}. \) Both boundary components of $A$ are invariant, and therefore $A$ is positively invariant. Let \( x\in A\setminus W^s(CE_+). \) Its positive semiorbit is contained in the bounded set $A$, and hence $\omega(x)$ is nonempty, compact, connected, and invariant. The set $\omega(x)$ cannot contain $CE_-$, because $CE_-\in D_{\rm hom}$, and it cannot contain $CE_+$ because $x\notin W^s(CE_+)$. Moreover, by hypothesis there is no periodic orbit in $A$ other than $\Gamma_H$. Since $\Gamma_{\rm hom}$ does not attract trajectories from $A$, the Poincar\'e-Bendixson theorem gives \( \omega(x)=\Gamma_H \) so that 

\[ 
\lim_{t\to+\infty} \operatorname{dist}\bigl(\Phi_t(x),\Gamma_H\bigr)=0 \qquad \forall x\in D_H\setminus \{ \overline{D_{\rm hom}} \cup W^s(CE_+) \}. 
\] 
Finally, assume that both $\overline{\mathcal B(\Gamma_H)}$ and $\overline{\mathcal B(HFE)}$ have nonempty intersections with $\mathbb R_+^2\setminus\overline{D_H}$. These intersections are open, disjoint, and positively invariant. If they coexist in the connected component of \( \mathbb R_+^2\setminus\overline{D_H} \) adjacent to $\Gamma_H$, their common relative boundary is nonempty. Choose a point \mbox{$z \in \partial \mathcal B_{R_+^2\setminus\overline{D_H}}(\Gamma_H) \cap \partial \mathcal B(HFE)$}. By boundedness of the positive solutions, $\omega(z)$ is a nonempty compact invariant set. Setting \( K:=\omega(z), \) one obtains \( K\subset \mathbb R_+^2\setminus D_H \) and \( K\subset \partial\mathcal B(\Gamma_H) \cap \partial\mathcal B(HFE). \) Moreover, $K$ is distinct from both $\Gamma_H$ and $HFE$, since an asymptotically stable invariant set belongs to the interior of its basin. Thus an additional compact invariant structure separates the asymptotic regimes governed by $\Gamma_H$ and $HFE$. \\

\noindent \textit{Proof of statement~{\rm (iv)}.}\\
Assume that \( D_{\rm hom}\subset D_H\) and \( \alpha>\alpha_c \) so that $\Gamma_{\rm hom}$ is repelling from the $D_{\rm hom}$ side and choose \( x\in D_{\rm hom}\setminus\{CE_-\}. \) The positive semiorbit of $x$ is bounded because $D_{\rm hom}$ is bounded and positively invariant. Therefore $\omega(x)$ is nonempty, compact, connected, and invariant. Since $CE_-$ is a hyperbolic repeller, \( \omega(x)\neq\{CE_-\}. \) Moreover, \( W^s(CE_+)\cap D_{\rm hom}=\varnothing, \) so $\omega(x)$ cannot contain $CE_+$. Since the homoclinic loop is repelling from the $D_{\rm hom}$ side, \( \omega(x)\neq\Gamma_{\rm hom}. \) It follows that $\omega(x)$ contains no equilibrium. By the Poincar\'e--Bendixson theorem, $\omega(x)$ is a periodic orbit contained in $D_{\rm hom}$. Define \( K:=\omega(x)\) a compact invariant set satisfying \( K\subset D_{\rm hom}\setminus\{CE_-\}. \) Since \( D_{\rm hom}\subset D_H\) and \( \Gamma_H=\partial D_H, \) one has \( \Gamma_H\cap D_{\rm hom}=\varnothing. \) Consequently, \( K\neq\Gamma_H. \) Thus $D_{\rm hom}$ contains at least one periodic orbit distinct from the Hopf cycle $\Gamma_H$.\\ 

\begin{figure}[h!] 
\centering 
\includegraphics[width=0.7\textwidth]{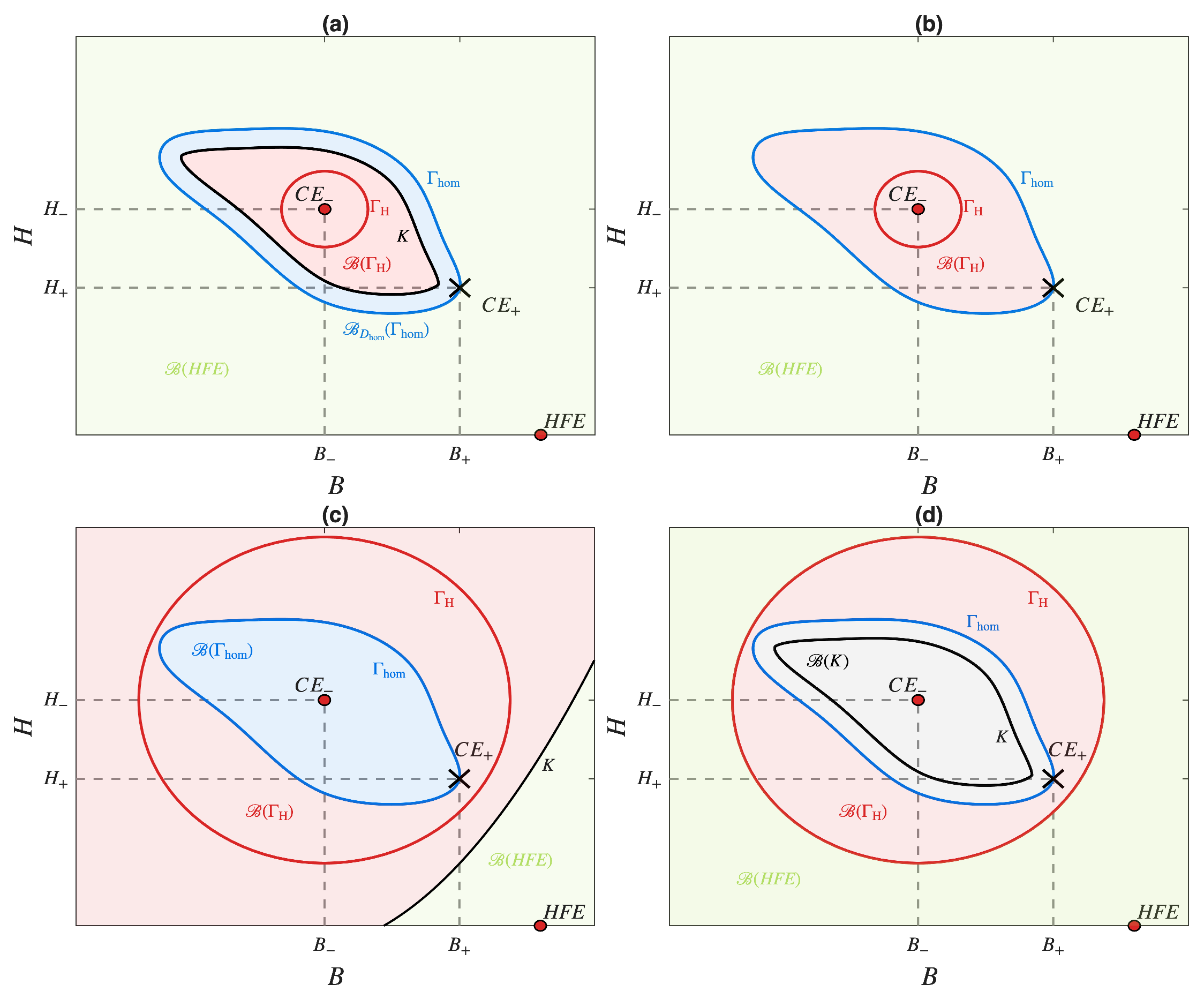} 
\caption{Schematic phase-plane sketches illustrating the qualitative organization of the invariant structures considered in Theorem~\ref{th:homoclinic_regionIII}. Panels (a)-(d) correspond to cases (i)-(iv), respectively.} 
\label{fig:hom_regionIII}
\end{figure}

A qualitative sketch of the above cases is provided in Figure~\ref{fig:hom_regionIII}. Finally, if $\alpha=\alpha_c$, then \( \sigma_+=0. \) In this case the leading-order term determining the stability of the homoclinic return map vanishes. Therefore the stability of $\Gamma_{\rm hom}$ cannot be inferred from the saddle quantity alone and requires a higher-order analysis of the Poincar\'e return map.
\end{proof}

\newpage 

\subsubsection{Section summary}

\begin{figure}[h!] 
\centering 
\includegraphics[width=0.95\textwidth]{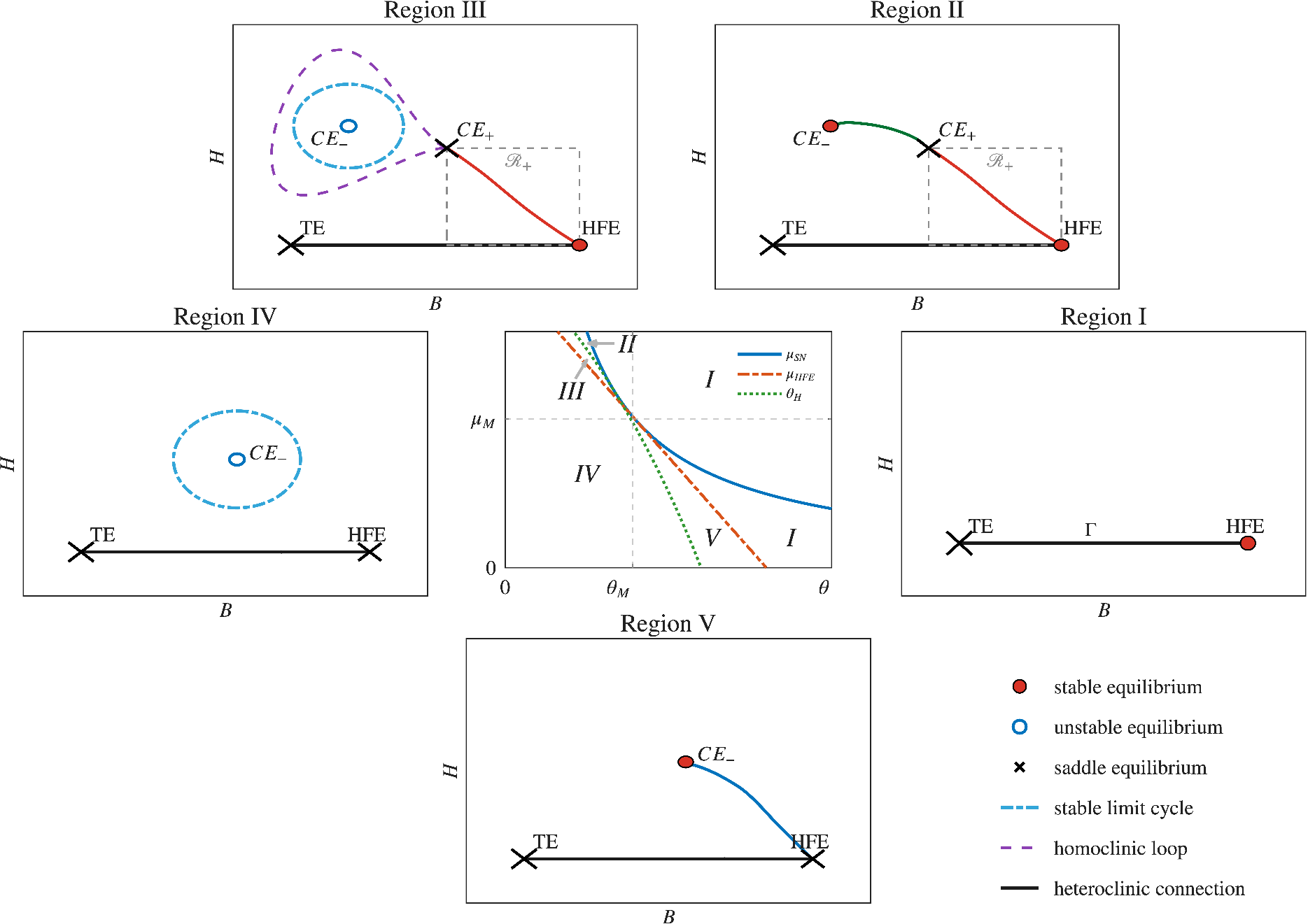} 
\caption{ Qualitative phase-plane configurations associated with the five regions of the parameter plane. Red dots denote stable equilibria, blue circles unstable equilibria, black crosses saddle points, dash-dotted blue curves stable limit cycles, dashed purple curves homoclinic loops, and solid curves heteroclinic connections.} 
\label{fig:global-dynamics-summary} 
\end{figure}

The classification established in Proposition~\ref{th:global}, together with the heteroclinic and homoclinic results of Theorems~\ref{thm:heteroclinics}-\ref{th:homoclinic_regionIII}, provides a complete qualitative description of the global phase-plane organization of system~\eqref{eq:model}. The corresponding dynamical configurations are summarized in Figure~\ref{fig:global-dynamics-summary}, where each panel represents the phase portrait associated with one of the five regions of the parameter plane. In particular:

\begin{enumerate}
\item[i)] Region~I corresponds to the simplest dynamical regime. The only biologically admissible equilibria are the trivial equilibrium $TE$ and the herbivore-free equilibrium $HFE$. The invariant boundary segment $\Gamma$ provides the heteroclinic connection $TE\to HFE$, while $HFE$ is globally attracting within the biologically relevant region int\,$\Omega$. Consequently, all trajectories eventually converge to herbivore extinction. 

\item[ii)] Region~II represents a rich equilibrium configuration in which both coexistence equilibria $CE_{-}$ and $CE_{+}$ are present, with $CE_{-}$ asymptotically stable and $CE_{+}$ a saddle. Besides the boundary connection $TE\to HFE$, the unstable manifold of $CE_{+}$ generates heteroclinic trajectories toward both $HFE$ and $CE_{-}$.

\item[iii)] Region III is characterized by the loss of stability of $CE_{-}$ through the supercritical Hopf bifurcation described in Theorem \ref{th:hopf}. As a consequence, a stable limit cycle $\Gamma_H$ surrounds the unstable coexistence equilibrium $CE_{-}$. The saddle equilibrium $CE_{+}$ persists and still generates the heteroclinic connection toward $HFE$. The global dynamics may, however, exhibit several additional configurations depending on the possible occurrence and geometric arrangement of homoclinic structures associated with $CE_{+}$. These scenarios, together with their dynamical implications, are analyzed in detail in Theorem \ref{th:homoclinic_regionIII} (or, in a qualitative sketch, in Figure \ref{fig:hom_regionIII}). For the sake of clarity, Figure \ref{fig:global-dynamics-summary} displays only a representative qualitative phase portrait of Region III rather than all admissible configurations. Consequently, the long-term dynamics are generally governed by stable periodic oscillations and/or alternative attractors separated by the invariant manifolds associated with $CE_+$, while additional invariant structures may arise in the presence of homoclinic configurations (see Theorem \ref{th:homoclinic_regionIII}).

\item[iv)] Region IV is characterized by the disappearance of the saddle equilibrium $CE_{+}$, while the coexistence equilibrium $CE_{-}$ remains unstable. As a consequence, the stable periodic orbit generated through the Hopf bifurcation persists and becomes the unique attractor in the biologically relevant region. The only heteroclinic connection supported by the system is the invariant boundary orbit $TE \to HFE$.

\item[v)] Region~V corresponds to a stable coexistence regime. The coexistence equilibrium $CE_{-}$ is asymptotically stable, whereas both $TE$ and $HFE$ are saddles. In addition to the boundary connection $TE\to HFE$, the unstable manifold of $HFE$ enters the interior of the phase plane and converges to $CE_{-}$, generating the heteroclinic connection $HFE\to CE_{-}$. Hence every biologically relevant trajectory belonging to int\,$\Omega$ converges to the coexistence equilibrium. 
\end{enumerate}

The five configurations illustrate how the interaction between local bifurcations and global invariant manifolds reorganizes the phase portrait. Depending on the parameter values, the system may exhibit herbivore extinction, stable coexistence, bistability, or persistent oscillatory dynamics, all of which are organized by a small number of heteroclinic, homoclinic, and periodic structures.

While Figure~\ref{fig:global-dynamics-summary} emphasizes the geometric organization of the phase portraits, Table~\ref{tab:global-dynamics} provides a synthetic overview of the invariant structures supported by the system in each parameter region. In particular, it highlights the attractors, the admissible heteroclinic connections, the possible homoclinic configurations, and the occurrence of stable limit cycles.

\begin{table}[h!] 
\centering 
\caption{Summary of the invariant structures and asymptotic dynamics in the five parameter regions.} 
\label{tab:global-dynamics} 
\renewcommand{\arraystretch}{1.3} \begin{tabular}{c c c c c } 
\hline 
\textbf{\small Region} & \textbf{\small Stable invariant sets} & \textbf{\small Heteroclinic con.} & \textbf{\small Homoclinic con.} & \textbf{\small Stable limit cycle} \\ 
\hline 
I & HFE & $TE\to HFE$ & None & None \\ 
\hline II & 
\begin{tabular}{c} 
$CE_{-}$\\ HFE 
\end{tabular} & 
\begin{tabular}{c} 
$TE\to HFE$\\ $CE_{+}\to HFE$\\ $CE_{+}\to CE_{-}$ 
\end{tabular} 
& None & None \\ 
\hline III & 
\begin{tabular}{c} 
$\Gamma_H$\\ 
$\Gamma_{\hom}$\\
HFE 
\end{tabular} & 
\begin{tabular}{c} 
$TE\to HFE$\\ $CE_{+}\to HFE$ 
\end{tabular} 
& $CE_{+}$ & Present \\ \hline IV & $\Gamma_H$ & $TE\to HFE$ & None & Present \\ \hline V & $CE_{-}$ & 
\begin{tabular}{c} 
$TE\to HFE$\\ $HFE\to CE_{-}$ 
\end{tabular} 
& None & None \\ 
\hline 
\end{tabular} 
\end{table}

\newpage 

\subsection{Numerical validation}

In this section, we provide a numerical validation of the analytical results established in Sections~\ref{sec:equilibria} and~\ref{sec:global}. Our objectives are twofold. First, we verify the local bifurcation scenario described by Proposition~\ref{prop:SN}, Theorem~\ref{th:hopf}, and Corollary~\ref{cor:hopf} through numerical continuation of equilibria and periodic solutions. Second, we reconstruct representative phase portraits associated with the five parameter regions identified in Proposition~\ref{th:global} and Theorems~\ref{thm:heteroclinics}-\ref{th:homoclinic_regionIII}, and compare the resulting invariant structures with the theoretical configurations summarized in Figure~\ref{fig:global-dynamics-summary}. 

Unless otherwise stated, the model parameters are selected so as to belong to the interior of the corresponding region in the parameter plane displayed in Figure~\ref{fig:phase_regions}. Direct numerical integrations of system~\eqref{eq:model} were performed in MATLAB using the adaptive Runge-Kutta solver \texttt{ode45} \cite{MATLAB}. Relative and absolute tolerances were fixed to \( \mathrm{RelTol}=10^{-8}\) and \( \mathrm{AbsTol}=10^{-10}, \) and all computations were repeated with stricter tolerances to verify numerical convergence. All simulations were carried out in the biologically admissible region \( \Omega=\mathbb{R}_{+}^{2}, \) and were integrated over sufficiently long time intervals to eliminate transient dynamics and identify the corresponding $\omega$-limit sets. Numerical bifurcation analysis was performed using the continuation package AUTO-07P \cite{Doedel1981}. Throughout the continuation procedure, the herbivore mortality parameter $\theta$ was used as the principal continuation parameter. 

\subsubsection{Numerical continuation of equilibria and periodic orbits} 
\label{sec:num_cont}

Let us begin by investigating the local bifurcation structure described in Section~\ref{sec:equilibria} through numerical continuation performed with AUTO-07P. The purpose of this analysis is to verify the existence and stability properties of the equilibrium branches predicted by Lemmas~\ref{lemma:steady}-\ref{lem:local-classification}, as well as the Hopf and saddle-node bifurcations established in Proposition~\ref{prop:SN} and Theorem~\ref{th:hopf}. On the other hand, we also corroborate the qualitative bifurcation scenario summarized in Figure~\ref{fig:qualitative-bifurcation-diagram}, showing that the sequence of stability transitions and bifurcation events predicted analytically is accurately reproduced by numerical continuation.

\begin{figure}[b!] 
\centering \includegraphics[width=1\textwidth] {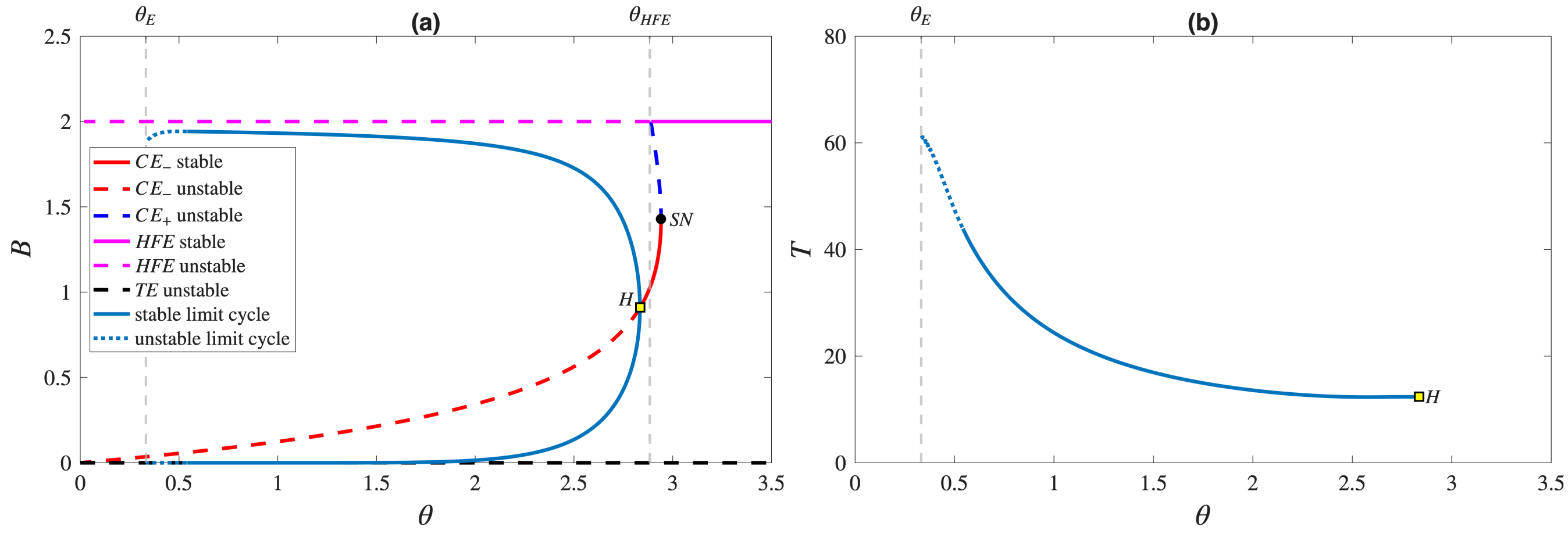} 
\caption{Numerical bifurcation diagram of system~\eqref{eq:model} with respect to the herbivore mortality parameter $\theta$, computed using AUTO--07P. (a) Equilibrium and periodic-solution branches. Solid and dashed curves denote stable and unstable branches, respectively. The coexistence equilibrium $CE_-$ loses stability at the Hopf point (H), generating a family of periodic orbits. The coexistence equilibria $CE_-$ and $CE_+$ collide at the saddle-node bifurcation (SN), while the herbivore-free equilibrium HFE changes stability at $\theta_{HFE}$. (b) Period $T$ of the periodic solutions as a function of $\theta$. Parameter values: $\alpha=2$, $m=10$ and $\mu=0.85$.} 
\label{fig:bifurcation-diagram} 
\end{figure}

In this context, Figure~\ref{fig:bifurcation-diagram} shows the bifurcation diagram obtained by using the herbivore mortality $\theta$ as continuation parameter. The diagram depicted in panel (a) reproduces all equilibrium branches predicted by the theoretical analysis. In particular, the herbivore-free equilibrium $HFE$ and the trivial state $TE$ persist for all parameter values, while the two coexistence equilibria $CE_{-}$ and $CE_{+}$ appear for intermediate values of $\theta$, namely for $0 < \theta \leq \theta_{SN}$ and $\theta_{HFE} < \theta \leq \theta_{SN}$, respectively. As expected from Lemma~\ref{lem:local-classification}, the branch $CE_{+}$ is always unstable, whereas the stability of $CE_{-}$ changes along the branch according to the position of the Hopf threshold. Indeed, the numerical continuation detects a Hopf bifurcation point (H) on the coexistence branch $CE_-$, where the equilibrium loses stability and a family of periodic solutions emerges. In agreement with Corollary~\ref{cor:hopf}, the bifurcating limit cycles are asymptotically stable for $\theta<\theta_H$ sufficiently close to the Hopf threshold. As $\theta$ is further decreased, the periodic branch remains stable over a wide parameter interval. However, the cycles are numerically detected as unstable along the terminal portion of the computed branch. Therefore, the supercritical nature of the Hopf bifurcation determines the local stability of the periodic solutions near $\theta_H$, but does not imply stability of the entire global family. The continuation terminates at $\theta=\theta_E$, where AUTO-07P reports an MX point. This label does not identify a bifurcation of the dynamical system, but indicates that the continuation could not be advanced within the prescribed numerical bounds. In particular, the endpoint may result from the increasing difficulty of resolving the periodic orbit as its period grows and its temporal profile becomes strongly nonuniform. Accordingly, $\theta_E$ should be interpreted as the terminal parameter value reached by the numerical continuation rather than as a bifurcation threshold. 

Panel (b) shows that the period $T$ of the oscillations increases rapidly as the branch approaches its terminal point. Such a growth is typically regarded as a classical signature of periodic solutions approaching a global bifurcation and suggests that the family may be converging towards one of the homoclinic configurations analysed in Section~\ref{subsec:hom}. This interpretation is supported by the geometric results obtained there. Indeed, Theorem~\ref{thm:no-homoclinic} identifies the saddle equilibrium $CE_+$ as the only equilibrium that can support a homoclinic loop, while Theorem~\ref{th:homoclinic_regionIII} characterizes the dynamical role of the resulting homoclinic structures in Region III. Therefore, the observed growth of the period, together with the numerical termination of the branch at the MX point, is consistent with the approach of the periodic orbit to a homoclinic configuration involving $CE_+$. Since neither the MX label nor the period growth alone provide a rigorous detection criterion, this interpretation should be viewed as numerical evidence rather than a proof of the existence of a homoclinic connection. The diagram also confirms the saddle-node bifurcation predicted by Proposition~\ref{prop:SN}. Indeed, the coexistence equilibria $CE_{-}$ and $CE_{+}$ approach each other and match at a turning point corresponding to the threshold $\theta_{SN}$. For larger values of $\theta$, no positive coexistence equilibrium exists and HFE becomes the only biologically admissible equilibrium. A further feature visible in Figure~\ref{fig:bifurcation-diagram} is the stability change of HFE at the threshold $\theta_{HFE}$. For $\theta<\theta_{HFE}$, HFE is unstable, whereas for $\theta>\theta_{HFE}$ it becomes locally asymptotically stable. Consequently, the interval $\theta_{HFE}<\theta<\theta_{SN}$ corresponds to a bistable regime in which HFE and $CE_{-}$ coexist as alternative attractors and the asymptotic outcome depends on the initial conditions. 

Overall, the numerical continuation reproduces all local bifurcation mechanisms derived analytically in Section~\ref{sec:equilibria} and provides quantitative confirmation of the theoretical bifurcation scenario underlying the global phase-plane organization investigated in the following subsections.

\subsubsection{Representative phase portraits in Regions I - V} 

\begin{figure}[t] 
\centering \includegraphics[width=1\textwidth] {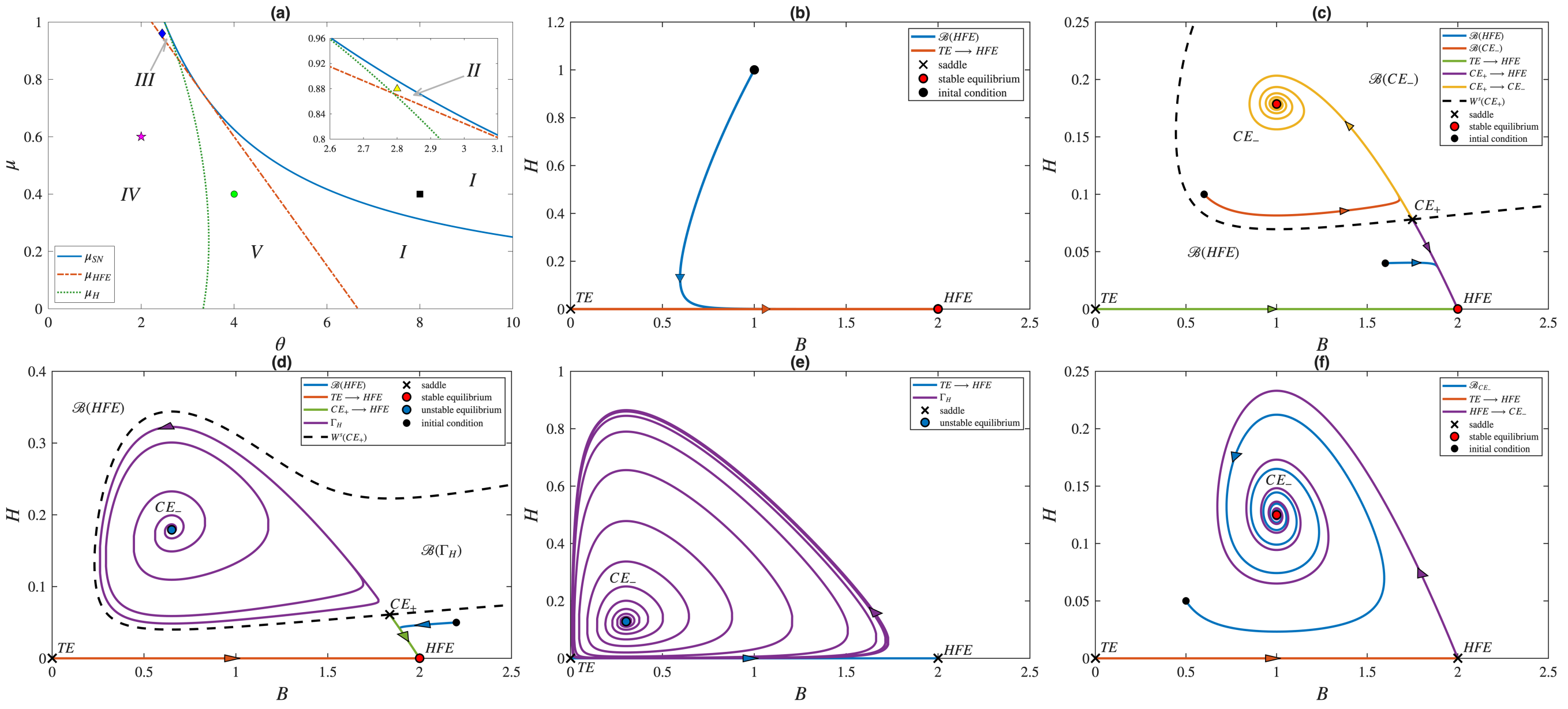} 
\caption{Numerical phase portraits associated with the five dynamical regions identified in Proposition \ref{th:global}. Panel (a) shows the partition of the parameter plane $(\theta,\mu)$ into Regions I - V. The markers indicate the parameter pairs $(\theta, \mu)$ selected for the corresponding phase portraits displayed in panels (b)-(f), respectively. Specifically, panel (b) represents Region I, panel (c) Region II, panel (d) Region III, panel (e) Region IV, and panel (f) Region V. Other parameter values fixed as in Figure \ref{fig:bifurcation-diagram}.} 
\label{fig:validation_phase_plane} 
\end{figure}

To complement the bifurcation analysis presented above, we now reconstruct representative phase portraits corresponding to the five parameter regions identified in Proposition \ref{th:global}. The aim is to verify numerically the qualitative phase-plane organization predicted by the theoretical analysis and sketched in Figure \ref{fig:global-dynamics-summary}, and to illustrate the invariant structures characterizing each dynamical regime. Figure \ref{fig:validation_phase_plane} summarizes these results by showing the partition of the parameter plane together with representative trajectories computed for parameter values selected within each region.

In particular, Figure \ref{fig:validation_phase_plane}(a) displays the partition of the parameter plane $(\theta,\mu)$ introduced in Proposition \ref{th:global}. The coloured markers identify the specific parameter pairs selected for the numerical simulations reported in panels (b)-(f). In particular, the black square, yellow triangle, blue diamond, magenta star, and green circle correspond to representative parameter values taken in Regions I - V, respectively. The inset provides a zoomed-in view of the parameter plane, where Region II becomes narrowly separated and difficult to distinguish at the scale of the full diagram.

Figure \ref{fig:validation_phase_plane}(b)-(f) display representative phase portraits corresponding to the five parameter regions identified in Proposition \ref{th:global}. Throughout these panels, black crosses denote saddles, red circles asymptotically stable equilibria, and blue circles unstable one. In particular, panel (b) corresponds to Region I in which, in agreement with Proposition \ref{th:global}, the only biologically admissible equilibria are the trivial equilibrium TE and the herbivore-free equilibrium HFE. The former is a saddle, whereas the latter is a stable node and represents the unique attractor of the system. The red trajectory reproduces the heteroclinic connection $TE \rightarrow HFE$ predicted by Theorem \ref{thm:heteroclinics}$_{(iii)}$, whereas the blue one, started from a generic point in the biologically relevant region, converges asymptotically to HFE, confirming that the basin of attraction of HFE coincides with the whole phase plane. Consequently, every biologically relevant solution approaches the herbivore-free equilibrium as $t \to +\infty$, validating the extinction scenario predicted by the theoretical analysis.

A markedly richer dynamical configuration is observed in panel (c), corresponding to Region II. In agreement with Proposition \ref{th:global}, the phase plane contains two saddle equilibria, namely TE and $CE_+$, together with the two asymptotically stable equilibria HFE and $CE_-$. The trajectories shown in the figure highlight both the heteroclinic connections predicted by Theorem \ref{thm:heteroclinics}$_{(iv)}$ and representative solution orbits converging towards the attractors of the system. In particular, the boundary orbit \mbox{TE $\rightarrow$ HFE} is recovered along the invariant axis $H=0$, while the unstable manifold of $CE_+$ generates the heteroclinic trajectories $CE_+$ $\rightarrow$ HFE and $CE_+$ $\rightarrow$ $CE_-$. The latter connection is consistent with the theoretical prediction in the absence of homoclinic loops associated with $CE_+$. The blue and red trajectories illustrate the bistable nature of this regime: depending on the initial condition, solutions converge either to HFE or to $CE_-$, confirming that both equilibria act as attractors. The stable manifold of the saddle equilibrium $CE_+$, displayed as a dashed black curve, clearly acts as the separatrix between the attraction basins of HFE and $CE_-$. Note that the stable invariant manifolds of the saddle equilibrium were reconstructed numerically from the eigenspaces of the linearized system starting from \( CE_+\pm \varepsilon v_u\) and \( CE_+\pm \varepsilon v_s, \) with $0 < \varepsilon \ll 1$, being $v_s$ and $v_u$ the normalized eigenvectors associated with the stable and unstable eigenvalues, respectively. Initial conditions located on opposite sides of this manifold converge to different attractors, in complete agreement with the analytical prediction of Theorem~\ref{thm:heteroclinics}.

A further qualitative transition is observed in panel (d), corresponding to Region III. In agreement with Proposition \ref{th:global}, the equilibria TE and $CE_+$ remain saddle points, HFE is still asymptotically stable, whereas the coexistence equilibrium $CE_-$ has lost stability through the Hopf bifurcation described in Theorem \ref{th:hopf}. As a consequence, a stable limit cycle emerges around $CE_-$, in accordance with Theorem \ref{th:hopf} and Corollary \ref{cor:hopf}. The phase portrait clearly reproduces the heteroclinic connections predicted by Theorem \ref{thm:heteroclinics}. In particular, the boundary orbit $TE \rightarrow HFE$ persists, while one branch of the unstable manifold of $CE_+$ generates the heteroclinic connection $CE_+ \rightarrow HFE$. At the same time, trajectories initialized in a neighbourhood of the unstable equilibrium $CE_-$ are repelled from it and converge towards the stable periodic orbit, confirming its attracting character. The numerical simulations therefore validate the theoretical prediction that Region III is dominated by persistent oscillatory dynamics. Consequently, the long-term dynamics are governed either by the stable limit cycle or by the herbivore-free equilibrium HFE, while the stable manifold of $CE_+$ acts as the separatrix between their attraction basins.

Panel (e) illustrates the subsequent dynamical regime arising after the disappearance of the saddle equilibrium $CE_+$. In agreement with Proposition \ref{th:global}, the phase plane now contains only the boundary equilibria TE and HFE, both of saddle type, together with the unstable coexistence equilibrium $CE_-$. The stable periodic orbit generated through the Hopf bifurcation in Region III persists after the disappearance of $CE_+$ and surrounds $CE_-$. 
In particular, it remains the unique attractor of the biologically relevant region. The numerical trajectories confirm this prediction. Solutions initiated from different points of the phase plane are first repelled from $CE_-$ and subsequently approach the stable limit cycle. Moreover, the only heteroclinic connection supported by the system is the boundary orbit $TE \rightarrow HFE$, consistently with Theorem \ref{thm:heteroclinics}. Consequently, the long-term dynamics are entirely dominated by persistent oscillations, validating the theoretical phase-plane organization predicted for Region IV.

Finally, panel (f) corresponds to Region V, where stable coexistence is recovered. In agreement with Proposition \ref{th:global}, the phase plane contains the two saddle equilibria TE and HFE together with the asymptotically stable coexistence equilibrium $CE_-$. The trajectories displayed in the phase portrait reproduce the heteroclinic structure predicted by Theorem \ref{thm:heteroclinics}. In particular, the boundary connection $TE \rightarrow HFE$ persists, while the unstable manifold of HFE enters the interior of the phase plane and generates the heteroclinic connection $HFE \rightarrow CE_-$. Furthermore, representative trajectories initialized in the biologically relevant region converge asymptotically to $CE_-$, confirming that the coexistence equilibrium constitutes the unique attractor of the system. Consequently, every biologically meaningful solution approaches a stable coexistence state between vegetation and herbivores, validating the theoretical phase-plane organization predicted for Region V.

\subsubsection{Numerical detection of homoclinic connections} 

The analytical results of Section~\ref{subsec:hom} identify the saddle equilibrium $\mathrm{CE}_{+}$ as the only equilibrium that may support a homoclinic connection and characterize the geometric configurations compatible with such a connection. The purpose of this section is to provide a numerical validation of these predictions. To this aim, we reconstruct the stable and unstable manifolds of $\mathrm{CE}_{+}$ and introduce a numerical splitting function measuring their relative position on a suitable transverse section. This approach allows us to illustrate the exclusion mechanism established in Theorem~\ref{prop:no-homoclinic-region-II} within Region~II and to detect a homoclinic connection in Region~III through a sign change of the splitting function, consistently with the criterion provided by Theorem~\ref{thm:homoclinic_criterion}.

To this aim, we introduce a transverse section 

\[ 
\Sigma= \{(B,H)\in\Omega:\ B=B_\Sigma\}, \qquad B_\Sigma<B_{+}, 
\] 
and follow the relative position of $W^{u}_{NW}(\mathrm{CE}_{+})$ and $W^{s}(\mathrm{CE}_{+})$ that, in this case, reduces to $W^{s}_{SW}(\mathrm{CE}_{+})$. Denoting by 

\[ 
q^{u}(\theta) = (B_\Sigma,H_\Sigma^{u}(\theta)) 
\qquad 
\text{and} 
\qquad
q^{s}(\theta) = (B_\Sigma,H_\Sigma^{s}(\theta)) 
\] 
their corresponding intersections with $\Sigma$, we define the numerical splitting function 

\[ 
\Delta_{\mathrm{num}}(\theta) = H_\Sigma^{u}(\theta) - H_\Sigma^{s}(\theta). 
\] 

\begin{figure}[t!] 
\centering \includegraphics[width=1\textwidth] {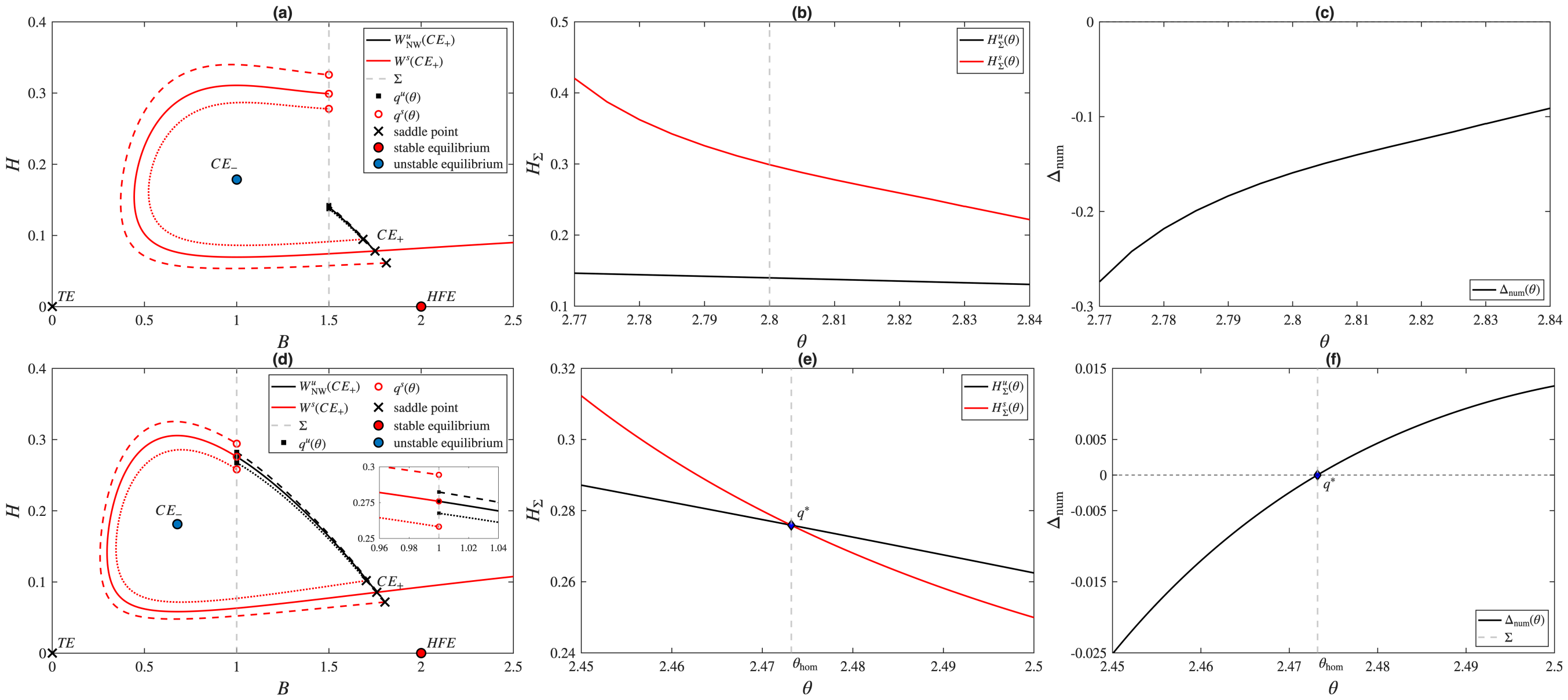} 
\caption{Numerical validation of the homoclinic framework of Section~\ref{subsec:hom}.The top row corresponds to a continuation path in Region~II, whereas the bottom one corresponds to Region~III. The left panels show the reconstructed stable manifolds $W^{s}(\mathrm{CE}_{+})$ (red) and unstable manifolds $W^{u}(\mathrm{CE}_{+})$ (black), together with the transverse section $\Sigma$ ($B_{\Sigma} = 1.5$ in panels (a)-(c) and $B_{\Sigma} = 1$ in panels (d)-(f)). In Region~II (III) the manifolds are displayed for $\theta=2.79$ ($\theta=2.46$) in dashed lines, $\theta=2.80$ ($\theta=2.4732$) in solid lines, and $\theta=2.81$ ($\theta=2.49$) in dotted one. The middle panels report the corresponding intersection coordinates $H_{\Sigma}^{u}$ and $H_{\Sigma}^{s}$ on the transverse section as functions of $\theta$. The right panels show the associated numerical splitting function $\Delta_{\mathrm{num}}(\theta)$.} 
\label{fig:validation_homoclinic} 
\end{figure}

According to Theorem~\ref{thm:homoclinic_criterion}, a zero of $\Delta_{\mathrm{num}}$ corresponds to an intersection between the selected stable and unstable manifold branches on the transverse section and therefore provides evidence of a homoclinic connection associated with $\mathrm{CE}_{+}$. Figure~\ref{fig:validation_homoclinic} summarizes the numerical results. The upper row corresponds to a continuation path entirely contained in Region~II (panels (a)-(c)), whereas the lower row corresponds to a continuation path in Region~III (panels (d)-(f)). Indeed, the first row provides a numerical illustration of the exclusion mechanism established analytically in Theorem~\ref{prop:no-homoclinic-region-II}. In particular, the reconstructed stable and unstable manifolds remain clearly separated throughout the biologically admissible continuation interval (see panel (a)) and the same behaviour is observed on the transverse section, where the intersection coordinates $H_\Sigma^{u}$ and $H_\Sigma^{s}$ remain distinct for all considered parameter values (see panel (b)). Consequently, the numerical splitting function $\Delta_{\rm num}$ does not approach zero and preserves a constant sign along the entire continuation path (panel (c)). Therefore, the selected manifold branches never exchange their relative position and no homoclinic configuration is detected. Although this calculation does not constitute a proof of nonexistence, it is fully consistent with the global analytical exclusion result provided by Theorem~\ref{prop:no-homoclinic-region-II}. 

On the other hand, the second row illustrates a qualitatively different situation arising in Region~III. As the parameter $\theta$ varies, the selected manifold branches approach each other and eventually exchange their relative position (see panel (d)). This transition is reflected in the section coordinates $H_\Sigma^{u}$ and $H_\Sigma^{s}$, which intersect at a parameter value close to \( \theta_{\mathrm{hom}} \approx 2.4732  \) (see panel (e)). Correspondingly, the splitting function $\Delta_{\rm num}$ changes sign and vanishes at an intermediate parameter value. The observed sign reversal indicates that the unstable and stable manifold branches pass through each other on the transverse section (see blue diamond in panels (e) and (f)), exactly as predicted by the geometric criterion of Theorem~\ref{thm:homoclinic_criterion}. The manifold reconstruction at $\theta\approx\theta_{\mathrm{hom}}$ simultaneously reveals the formation of a homoclinic loop based at $\mathrm{CE}_{+}$. These observations provide numerical evidence of a codimension-one homoclinic connection in Region~III. 

The numerical computations of Figure~\ref{fig:validation_homoclinic} therefore reproduce the two fundamental geometric scenarios predicted by the theoretical analysis: (i) in Region~II the manifold branches remain separated and the splitting function stays away from zero, consistently with Theorem~\ref{prop:no-homoclinic-region-II}, so that no homoclinic solution emerges; (ii) in Region~III the splitting function undergoes a sign change and a homoclinic loop is observed, in agreement with the criterion established in Theorem~\ref{thm:homoclinic_criterion}. 

\section{Conclusion}
\label{sec:disc}

In this work, we investigated the global dynamics of a toxin-mediated plant-herbivore model with a non-monotone defence-dependent functional response. 

Building upon the local bifurcation analysis previously developed for this class of models, we focused on the geometric structures that organize the phase portrait and determine the long-term behaviour of the ecological system. First, we provided a complete characterization of the biologically admissible equilibria and their stability properties. The analysis revealed how the transition from a monotone to a non-monotone consumption function generates qualitatively distinct equilibrium configurations and partitions the parameter space into different dynamical regimes. In particular, the coexistence equilibrium $CE_{-}$ may lose stability through a supercritical Hopf bifurcation, leading to the emergence of stable periodic oscillations. The main contribution of the paper is the complete geometric classification of the global phase-plane organization of the model, including the characterization of the heteroclinic framework, the admissible homoclinic configurations, and the associated parameter-space partition. In particular, we obtain a complete classification of the heteroclinic connections supported by the model and derive rigorous constraints on the occurrence of homoclinic configurations. These results complement the existing local bifurcation theory by revealing how equilibrium states, periodic orbits, and global connections are embedded within a unified phase-plane structure. From an ecological perspective, the analysis clarifies how toxin-mediated inhibition can produce markedly different outcomes, including herbivore extinction, stable coexistence, and persistent population oscillations. More importantly, it shows that these ecological regimes are organized by a relatively small set of invariant structures whose configuration changes as the model parameters vary. 

Several questions remain open. The most immediate one concerns the multiscale nature of the oscillatory regime. Numerical evidence suggests that some periodic solutions may display slow-fast behaviour, and a rigorous analysis within the framework of Geometric Singular Perturbation Theory could provide further insight into the mechanisms underlying these oscillations \cite{Iuorio2022, Iuorio2023, Grifo2025IV, Borsotti2025}. Then, the interplay between spatial effects and multiscale dynamics may generate complex spatiotemporal behaviours, including relaxation oscillations, localized structures, and pattern transitions occurring on multiple temporal and spatial scales. Developing a unified geometric framework for these phenomena represents a promising direction for future research \cite{Grifo2025, Grifo2025II, Grifo2025III, Grifo2026}.

\section*{Acknowledgments}
This research was funded by INdAM-GNFM.

\section*{Data availability}
No data was used for the research described in the article.

\appendix 

\section{Sign of the first Lyapunov coefficient} 
\label{app:first_lyapunov}
The purpose of this appendix is to prove that the first Lyapunov coefficient associated with the Hopf bifurcation occurring along the coexistence branch $\mathrm{CE}_{-}$ is always negative. Recall that, after imposing the Hopf condition, the quantity determining the sign of the first Lyapunov coefficient can be written as 

\[ 
\Delta_H = \frac{4B_H} {(1+B_H)^7 \left[(1-\mu)B_H+\mu+1\right]} \,P(B_H,\mu), 
\] 
where 

\[ 
P(B_H,\mu) =(1-\mu)^2(2\mu-1)B_H^4 +(1-\mu)(\mu^2+6\mu-4) B_H^3+3(\mu^2+2\mu-2)B_H^2 +(3\mu^2-2\mu-4)B_H -(2\mu+1). 
\] 
Since $B_H>0$ and $1/4 < \mu < 1$, the sign of $\Delta_H$ coincides with the sign of $P$. 

Let us first assume that $1/4<\mu\le1/2$. Since all the coefficients here involved are non positive, it is easy to ascertain that $P(B_H,\mu) < 0$ for any $\mu \in [1/4, 1/2]$. 

Then, to prove the assert, let us now assume that $1/2 < \mu \leq 1$. Since the equilibrium $\mathrm{CE}_{-}$ belongs to the increasing branch of $C(B)$, one has $ C'(B_H)>0 $ so that $0<(2\mu-1)B_H<1$. By introducing the variables $s:=2\mu-1$ and  $t:=sB_H$, it follows that $0<s\le1$, $0<t<1$, $\mu=(1+s)/2$, and $B_H=t/s$. Then, defining 

\[ 
F(s,t) := 16s^4 P\!\left( \frac{t}{s}, \frac{1+s}{2} \right),
\] 
it can be proved that, since $s>0$, the sign of $P(B_H,\mu)$ is entirely determined by the sign of $F(s,t)$. Therefore, direct computations yield 

\[ 
F(s,t) = 4s(1-s)^2 t^4 - 2 s \bigl( s^3+13s^2-17s+3 \bigr)t^3 +12 s^2 \bigl(s^2+6s-3 \bigr)t^2 + 4 s^3 \bigl( 3s^2+2s-17 \bigr)t -16 s^4 \left( s + 2 \right) . 
\] 
Since $F$ has degree five in $s$ and degree four in $t$, it admits the tensor-product Bernstein representation \cite{lorentz1986}

\[ 
F(s,t) = \sum_{i=0}^{5} \sum_{j=0}^{4} b_{ij} T_i^5(s)\, T_j^4(t), 
\] 
where 

\[ 
T_i^5(s) = \binom{5}{i} s^i(1-s)^{5-i}, 
\qquad 
T_j^4(t) = \binom{4}{j} t^j(1-t)^{4-j}. 
\] 
Therefore, the Bernstein coefficients are denoted by 

\[ 
[b_{ij}] = 
\begin{bmatrix} 
0&0&0&0&0\\[1mm] 
0&0&0&-\frac{3}{10}&-\frac25\\[1mm] 
0&0&-\frac35&-\frac{31}{20}&-\frac95\\[1mm] 
0&-\frac{17}{10}&-4&-\frac{59}{10}&-6\\[1mm] 
-\frac{32}{5}&-\frac{64}{5}&-\frac{88}{5}&-\frac{98}{5}&-\frac{88}{5}\\[1mm] 
-48&-60&-64&-60&-48 
\end{bmatrix}. 
\] 

Since the Bernstein bases $T_i^5(s)>0$ and $T_j^4(t)>0$ for every $(s,t)\in]0,1[ \times ]0,1[$, the sign of $F(s,t)$ is entirely determined by the Bernstein coefficients. In our case, they satisfy $b_{ij}\le0$ where some of them are strictly negative, e.g. $b_{40}=-32/5<0$, so that it follows that $F(s,t)<0$ on $]0,1] \times ]0,1[$  and, consequently, $P(B_H,\mu)<0$, leading to $l_1(\theta_H)<0$.

\bibliographystyle{unsrtnat}
\bibliography{bibliography}

\end{document}